\documentclass[11pt]{article}
\usepackage{amsfonts}
\usepackage{mathrsfs}
\usepackage{amssymb}
\usepackage{amsmath}
\usepackage{amsthm}
\usepackage{latexsym}
\usepackage{graphicx}
\usepackage{color}
\usepackage{hyperref}
\usepackage{ulem}
\usepackage{fancyhdr}
\usepackage{lastpage}
 \numberwithin{equation}{section}

\newcommand{\R}{\mathbb{R}^3\,}

\newcommand{\p}{\partial}
\newcommand{\e}{\varepsilon}
\newcommand{\dd}{{\rm d}}

\newtheorem{Theorem}{Theorem}[section]
\newtheorem{Lemma}{Lemma}[section]

\newtheorem{proposition}{Proposition}[section]
\newtheorem{definition}{Definition}[section]

\allowdisplaybreaks[4]

\begin{document}

\begin{center}

        \vspace{1cm}
        \LARGE{\textbf{A Kato-Type Criterion for Vanishing Viscosity and Magnetic Diffusivity Limits for Compressible Magnetohydrodynamics}}
        \vspace{1cm}

        \large
        {\bf Siran Li}\\
        School of Mathematical Sciences $\&$ CMA-Shanghai\\
        Shanghai Jiao Tong University, Shanghai 200240, China\\
        E - mail: siran.li@sjtu.edu.cn\\
        \vspace{1cm}
        {\bf Tong Tang}\\
        Department of Mathematics, College of Sciences\\
        Yangzhou University, Yangzhou 225000, China\\
        E - mail: tt0507010156@126.com\\
         \vspace{1cm}
       {\bf Yuanjing Tang}\\
        Department of Mathematics, College of Sciences\\
        Yangzhou University, Yangzhou 225000, China\\
        E - mail: Tyj18118851855@163.com\\
        \vspace{0cm}

        \rule{\textwidth}{1pt}
    \end{center}

\begin{abstract}
In this paper, we investigate the vanishing viscosity and magnetic diffusivity limits for the compressible magnetohydrodynamics (MHD) equations on a smooth bounded domain in $\mathbb{R}^3$. Under the assumptions of interior regularity and uniform local boundedness near the boundary, we establish a Kato-type criterion under which finite energy weak solutions to the dissipative MHD system converge to weak solutions to the ideal MHD system as the viscosity and magnetic resistivity coefficients decrease to zero. In particular, under such criterion, the anomalous energy dissipation in the vanishing viscosity and diffusivity limits for MHD is precluded.
\end{abstract}

\textbf{Mathematics Subject Classification:} 35Q85, 76N10, 76N20, 76W05

\vspace{0.2 cm}

\textbf{Keywords:} inviscid and non-diffusive limit; compressible magnetohydrodynamics; Kato-type criterion; anomalous dissipation.

\section{Introduction}

Let $\Omega \subset \R$ be a smooth bounded domain. The partial differential equations for magnetohydrodynamics (MHD), first introduced by Alfv\'{e}n in the seminal work~\cite{al}, are fundamental equations governing the behavior of electrically
conducting fluids on $\Omega$ in the presence of magnetic fields:
\begin{equation}
\label{1.1}
\begin{cases}
\partial_t \rho + \operatorname{div}(\rho \mathbf{u}) = 0,\\
\partial_t (\rho \mathbf{u}) + \operatorname{div}(\rho \mathbf{u} \otimes \mathbf{u}) + \nabla P = \operatorname{curl} \mathbf{B} \times \mathbf{B} + \nu\Delta\mathbf u,\\
\partial_t \mathbf{B} - \operatorname{curl}(\mathbf{u} \times \mathbf{B}) = - \operatorname{curl}(\eta\operatorname{curl} \mathbf{B}),\\
\operatorname{div} \mathbf{B} = 0 \qquad\qquad\text{ in } [0,T] \times \Omega;  \\
(\rho, \mathbf{u}, \mathbf{B})|_{t=0} = (\rho_0, \mathbf{u}_0, \mathbf{B}_0)\qquad\qquad\text{at } \{0\} \times \Omega.
\end{cases}
\end{equation}
Here, $\rho$ and $\mathbf{u}$ denote the density and velocity of the compressible fluid, and $\mathbf{B}$ represents magnetic field experienced by the fluid. The Lorentz force exerted on the fluid is $$F_{\rm Lorentz} = \operatorname{curl}\mathbf{B} \times \mathbf{B}.$$ For an isentropic gas, the pressure $P$ satisfies
\begin{align*}
P(\rho) = \rho^\gamma
\end{align*}
with the adiabatic exponent \(\gamma > 1\). The  viscosity and magnetic diffusivity are $\nu$ and $\eta$, respectively.

The compressible MHD system is derived from the first principles of continuum mechanics and electromagnetism: the velocity field obeys the Navier--Stokes equations with the Lorentz force, and the magnetic field satisfies the Maxwell--Faraday equations. Unlike its incompressible counterpart, Eq.~\eqref{1.1}  captures the coupling between thermodynamic pressure and magnetic forces, hence is suitable for describing important wave propagation phenomena including fast/slow magnetosonic waves, shocks in fusion plasmas, as well as stellar atmospheres and interstellar turbulence, etc.

In typical settings of astronomy and plasma flows~\cite{du,sh}, the MHD system is often characterized by large hydrodynamic and magnetic Reynolds numbers, or equivalently, small viscosity and magnetic resistive coefficients. This motivates our mathematical investigation in the regime $$\nu=\eta=\varepsilon.$$ This recasts Eq.~\eqref{1.1} into the following:
\begin{equation}
\label{1.2}
\begin{cases}
\partial_t \rho^\varepsilon + \operatorname{div}(\rho^\varepsilon \mathbf{u}^\varepsilon) = 0,\\[4pt]
\partial_t (\rho^\varepsilon \mathbf{u}^\varepsilon) + \operatorname{div}(\rho^\varepsilon \mathbf{u}^\varepsilon \otimes \mathbf{u}^\varepsilon) + \nabla (\rho^\varepsilon)^\gamma = \operatorname{curl} \mathbf{B}^\varepsilon \times \mathbf{B}^\varepsilon + \varepsilon \Delta \mathbf{u^\varepsilon},\\[4pt]
\partial_t \mathbf{B}^\varepsilon + \varepsilon \operatorname{curl}(\operatorname{curl} \mathbf{B}^\varepsilon) = \operatorname{curl}(\mathbf{u}^\varepsilon \times \mathbf{B}^\varepsilon),\\[4pt]
\operatorname{div} \mathbf{B}^\varepsilon = 0 \qquad\qquad\text{ in } [0,T] \times \Omega;  \\
(\rho^\varepsilon, \mathbf{u}^\varepsilon, \mathbf{B}^\varepsilon)|_{t=0} = (\rho_0, \mathbf{u}_0, \mathbf{B}_0)\qquad\qquad\text{at } \{0\} \times \Omega,
\end{cases}
\end{equation}
supplemented with the boundary conditions:
\[
\mathbf{u}^\varepsilon|_{\partial\Omega} = \mathbf{0}, \quad
\mathbf{B}^\varepsilon \cdot \mathbf{n}|_{\partial\Omega} = 0, \quad
\operatorname{curl} \mathbf{B}^\varepsilon \times \mathbf{n}|_{\partial\Omega} = \mathbf{0}.
\]
Here and hereafter, \(\mathbf{n}\) denotes the unit outward normal vector field along $\p\Omega$.

In physical experiments, numerical analysis, and mathematical fluid dynamics, the vanishing viscosity limit (\textit{a.k.a.} the inviscid limit)  is a central topic of research. For the Navier--Stokes equations, namely by setting ${\bf B}^\e \equiv {\bf 0}$ in \eqref{1.2}, the inviscid limit is a highly non-trivial singular perturbation problem, primarily due to the mismatch between the no-slip boundary condition for the Navier--Stokes equations and the ideal non-penetration condition for the limiting system, \textit{i.e.}, the Euler equations. Such discrepancy generates sharp Prandtl-type boundary layers that may induce anomalous energy dissipation, which is crucial for the study of turbulence. A successful  theory developed by Kolmogorov in 1941 for the incompressible Navier--Stokes equations, known as the K41 theory, asserts among other results that the dimensionless dissipation rate approaches a positive constant in the inviscid limit.

Onsager~\cite{o} conjectured that there exist weak solutions to the incompressible Euler system of regularity $C^{0,\alpha}$ that are non-unique and dissipate kinetic energy, provided that $\alpha < 1/3$. Based on pioneering works by De Lellis--Sz\'{e}kelyhidi~\cite{de} and later developments (\textit{e.g.},  Buckmaster--De Lellis--Sz\'{e}kelyhidi--Vicol~\cite{bu}, among many others), which in turn derive motivations from convex integration pioneered by Nash in his groundbreaking work on $C^1$-isometric embeddings~\cite{nash}, Isett~\cite{i} finally settled Onsager's conjecture. In fact, it is known that there exists an abundance of ``wild''  weak solutions to the incompressible Euler equations. Thus, a prominent question is the \emph{selection principle} for physically reasonable weak solutions, of which a well-received candidate is the  vanishing viscosity limit.

Mathematical investigations on the inviscid limit are abundant in the literature; see, \textit{e.g.},  \cite{bardos19,ma,lo,sam,chr}. In this paper, we focus on an important type of criteria pioneered by Kato~\cite{ka} and apply it for the first time, to the best of our knowledge, to the compressible MHD system.

The Kato-type condition ascertains that if the energy dissipation tends to zero in a boundary layer of width comparable to the viscosity, then the inviscid limit is valid. Since its proposal by Kato~\cite{ka} in 1984, continuous efforts have been made towards its improvement and extension; see, \textit{e.g.}, \cite{wang,ke}.  For the compressible Navier--Stokes equations, Sueur~\cite{su} and Bardos--Nguyen~\cite{bardos16} established variants of Kato's criterion, and Wang--Zhu~\cite{wang1} improved these results by using only certain components of the velocity field. All of these works address the convergence of weak solutions to the Navier-Stokes equations to  strong solutions to the Euler equations. The main tools for the proofs include the arguments for the relative entropy and artificial boundary layer.  We refer the reader  to \cite{bi,fan,faraco20,f2,f3,liu}, among many other references, for relevant developments.

On the other hand, much less is known for the inviscid limit when considering weak solutions to the limiting Euler system. In the seminal work~\cite{cgq}, Chen--Perepelitsa proved that the solutions to the 1D compressible Navier--Stokes equations converge to an entropy weak solution of finite-energy to the Euler system. Chen--Liang--Wang~\cite{chen2022} established a Kato-type criterion for the incompressible Navier--Stokes equations (see also the earlier work by Drivas--Nguyen~\cite{dr}), by way of introducing a novel boundary layer foliation scheme. The same authors, together with Xu, extended in \cite{chen2024} the above result to the inviscid limit of the compressible Navier--Stokes equations.

Inspired by the aforementioned developments, in this paper we investigate the inviscid and non-resistive limits of the compressible MHD system. The limiting PDE, obtained by formally setting $\varepsilon = 0$ in \eqref{1.2}, is the ideal compressible MHD system: 
\begin{equation}
\label{1.3}
\begin{cases}
\partial_t \rho^0 + \operatorname{div}(\rho^0 \mathbf{u}^0) = 0,\\[4pt]
\partial_t (\rho^0 \mathbf{u}^0) + \operatorname{div}(\rho^0 \mathbf{u}^0 \otimes \mathbf{u}^0) + \nabla (\rho^0)^\gamma = (\operatorname{curl} \mathbf{B}^0) \times \mathbf{B}^0,\\[4pt]
\partial_t \mathbf{B}^0 = \operatorname{curl}(\mathbf{u}^0 \times \mathbf{B}^0),\\[4pt]
\operatorname{div} \mathbf{B}^0 = 0
\qquad\qquad\text{ in } [0,T] \times \Omega, \\[4pt]
(\rho^0, \mathbf{u}^0, \mathbf{B}^0)|_{t=0} = (\rho_0, \mathbf{u}_0, \mathbf{B}_0) \qquad\qquad\text{at } \{0\} \times \Omega,
\end{cases}
\end{equation}
supplemented by the boundary conditions
\[
\mathbf{u}^0 \cdot \mathbf{n}|_{\partial\Omega} = 0, \quad
\mathbf{B}^0 \cdot \mathbf{n}|_{\partial\Omega} = 0.
\]


Recall the definition of the renormalized weak solutions for the dissipative MHD system \eqref{1.2}, which is an important subclass of classical weak solutions. The existence theory of these solutions may be seen, \textit{e.g.}, in \cite{HuWang2010} and \cite[Chapter~3]{f1}.

\begin{definition}\label{def:weak_solution}
Given a bounded domain $\Omega \subset \mathbb{R}^3$ and an adiabatic exponent $\gamma > \frac{3}{2}$. For any fixed $\varepsilon > 0$ and initial data $(\rho^\varepsilon_0, \mathbf{m}^\varepsilon_0, \mathbf{B}^\varepsilon_0)$ satisfying $\rho^\varepsilon_0 \ge 0$, $\rho^\varepsilon_0 \in L^\gamma(\Omega)$, $\frac{|\mathbf{m}^\varepsilon_0|^2}{\rho^\varepsilon_0} \in L^1(\Omega)$ with $\mathbf{m}^\varepsilon_0 = \mathbf{0}$ a.e. on $\{\rho^\varepsilon_0 = 0\}$, and $\mathbf{B}^\varepsilon_0 \in L^2(\Omega)$ with $\operatorname{div}\mathbf{B}^\varepsilon_0 = 0$ in $\mathcal{D}'(\Omega)$, there exists a \emph{renormalized weak solution} $(\rho^\varepsilon, \mathbf{u}^\varepsilon, \mathbf{B}^\varepsilon)$ to the initial-boundary value problem \eqref{1.2}, such that:
\begin{itemize}
    \item The density $\rho^\varepsilon \ge 0$ satisfies
    $$
    \rho^\varepsilon \in C\left([0,T]; L^1(\Omega)\right) \cap L^\infty\big([0,T]; L^\gamma(\Omega)\big), \qquad \rho^\varepsilon(x,0) = \rho^\varepsilon_0,
    $$
    and the momentum $\rho^\varepsilon \mathbf{u}^\varepsilon$ satisfies
    $$
    \rho^\varepsilon \mathbf{u}^\varepsilon \in C\left([0,T]; L^{\frac{2\gamma}{\gamma+1}}_{\mathrm{weak}}(\Omega)\right);
    $$
    \item The velocity $\mathbf{u}^\varepsilon$ and the magnetic field $\mathbf{B}^\varepsilon$ satisfy
    \begin{align*}
    \mathbf{u}^\varepsilon \in L^2\left([0,T]; W_0^{1,2}(\Omega)\right), \quad \mathbf{B}^\varepsilon \in L^2\left([0,T]; W^{1,2}(\Omega)\right) \cap C_{\rm weak}\left([0,T]; L^2(\Omega)\right);
    \end{align*}
    \item The system \eqref{1.2} is satisfied in the sense of distributions, that is, in $\mathcal{D}'(\Omega \times (0,T))$;
        \item The continuity equation holds in the \emph{renormalized sense}: for any $b \in C^1(\mathbb{R}^+)$ satisfying
    $$
    b'(z) = 0 \quad \text{for all } z \ge z_0
    $$
    where $z_0>0$ is a constant depending only on  $b$, the renormalized equation
    $$
    \partial_t b(\rho^\varepsilon) + \operatorname{div}\bigl(b(\rho^\varepsilon)\mathbf{u}^\varepsilon\bigr) + \bigl(b'(\rho^\varepsilon)\rho^\varepsilon - b(\rho^\varepsilon)\bigr)\operatorname{div}\mathbf{u}^\varepsilon = 0
    $$
    holds in $\mathcal{D}'(\Omega \times (0,T))$;

    \item The energy inequality
    $$
    \mathcal{E}^\varepsilon(t) + \varepsilon \int_0^t \int_{\Omega} \left(|\nabla \mathbf{u}^\varepsilon|^2  + |\operatorname{curl} \mathbf{B}^\varepsilon|^2 \right) \dd x\,\dd s \le \mathcal{E}^\varepsilon(0)
    $$
    holds for a.e. $t \in [0, T]$, where
    $$
    \mathcal{E}^\varepsilon(t) = \int_{\Omega} \left( \frac{1}{2}\rho^\varepsilon |\mathbf{u}^\varepsilon|^2 + \frac{(\rho^\varepsilon)^\gamma}{\gamma-1} + \frac{1}{2}|\mathbf{B}^\varepsilon|^2 \right)(t,x)\, \dd x
    $$
    and
    $$
    \mathcal{E}^\varepsilon(0) = \int_{\Omega} \left( \frac{|\mathbf{m}^\varepsilon_0|^2}{2\rho^\varepsilon_0} + \frac{(\rho^\varepsilon_0)^\gamma}{\gamma-1} + \frac{1}{2}|\mathbf{B}^\varepsilon_0|^2 \right)(x)\, \dd x.
    $$
\end{itemize}
\end{definition}

Throughout this paper, for any $h>0$, we denote 
\begin{align*}
    \Omega^h := \{x \in \Omega : \operatorname{dist}(x, \partial\Omega) > h\}\quad\text{and}\quad \Gamma_h := \Omega \setminus \overline{\Omega^h}.
\end{align*}  
By ``uniform'' (boundedness/convergence/estimates) we mean independent of $\e$ as $\e\to 0$. Meanwhile, for simplicity, we shall refer to the vanishing viscosity and resistance limit $\e\to 0$ for the MHD system~\eqref{1.2} as the \emph{inviscid limit}.

Our main result in this paper, Theorem~\ref{thm: main}, provides a Kato-type criterion for the  inviscid limit of compressible MHD equations, assuming the spacetime Besov $B^{\sigma, \infty}_3$-regularity of the density $\rho^\e$, velocity ${\bf u}^\e$, and magnetic field ${\bf B}^\e$, as well as the uniform boundedness in $L^4_tL^\infty_x$ of $({\bf u}^\e, {\bf B}^\e)$ in the boundary layer $\Gamma_{8\e}$. One of the key ingredients of the proof is the domain foliation with variable width, introduced in Chen--Liang--Wang~\cite{chen2022} and Chen--Liang--Wang--Xu~\cite{chen2024} for the inviscid limit problem of the Navier--Stokes equations.

\begin{Theorem}\label{thm: main}
    Let $\Omega \subset \mathbb{R}^3$ be a bounded domain of $C^3$-regularity. Let $(\rho^\varepsilon, \mathbf{u}^\varepsilon, \mathbf{B}^\varepsilon)$ be a weak solution to the compressible MHD Equation~\eqref{1.2} such that
    \begin{equation}\label{1.4}
    \gamma \geqslant 2 \quad \text{and}\quad  0 \leqslant \rho^\varepsilon(x,t) \leqslant \widehat{\rho} < \infty
    \end{equation}
    for some constant $\widehat{\rho}>0$ independent of $\varepsilon$. Assume for any $\widetilde{\Omega}_T \subset\subset \Omega^\varepsilon \times (\varepsilon, T-\varepsilon)$ that
    \begin{align}\label{1.5}
    &\rho^\varepsilon, \ \mathbf{u}^\varepsilon \text{ are uniformly bounded in } B_3^{\sigma_1, \infty}(\widetilde{\Omega}_T) \text{ and } \nonumber\\
    &\mathbf{B}^\varepsilon \text{ is uniformly bounded in } B_3^{\sigma_2, \infty}(\widetilde{\Omega}_T), \text{ where }  \sigma_1, \sigma_2 \in \left(\frac{1}{3}, \frac{1}{2}\right],
    \end{align}
and for the boundary region $\Gamma_{8\varepsilon}$ that
    \begin{equation}\label{1.6}
        \mathbf{u}^\varepsilon \text{ and } \mathbf{B}^\varepsilon \text{ are uniformly  bounded in } L^4\big(0, T; L^\infty(\Gamma_{8\varepsilon})\big).
    \end{equation}
Then, the absence of anomalous dissipation, namely that 
    \begin{equation}\label{1.7}
        \lim_{\varepsilon \to 0} \varepsilon \int_0^T \int_\Omega \left( |\nabla \mathbf{u}^\varepsilon|^2 + |\nabla \mathbf{B}^\varepsilon|^2 \right) \,\dd x \,\dd t = 0,
    \end{equation}
holds whenever 
    \begin{equation}\label{1.8}
        \lim_{\varepsilon \to 0} \varepsilon \int_0^T \int_{\Gamma_{8\varepsilon}} \left( |\nabla \mathbf{u}^\varepsilon|^2 + |\nabla \mathbf{B}^\varepsilon|^2 \right) \,\dd x\,\dd t = 0.
    \end{equation}
    
    Moreover, after passing to subsequences, the solution $(\rho^\varepsilon, \mathbf{u}^\varepsilon, \mathbf{B}^\varepsilon)$ converges weakly to a weak admissible solution $(\rho^0, \mathbf{u}^0, \mathbf{B}^0)$ to the ideal compressible MHD Equation~\eqref{1.3} as $\varepsilon \to 0$. For any $p \in (1, \frac{9}{2}]$ one further has the strong convergence results:
    \begin{equation}\label{1.9}
        \rho^\varepsilon \longrightarrow \rho^0, \quad \rho^\varepsilon\mathbf{u}^\varepsilon \longrightarrow \rho^0\mathbf{u}^0, \quad \mathbf{B}^\varepsilon \longrightarrow \mathbf{B}^0 \quad \text{in } L^3\left(0, T; L^p_{\rm loc}(\Omega)\right).
    \end{equation}

\end{Theorem}

A few remarks are in order:
 
\begin{enumerate}
\item
As in Feireisl--Gwiazda--\'{S}wierczewska-Gwiazda--Wiedemann~\cite{f}, $\rho^\varepsilon,\mathbf u^\varepsilon\in  B_3^{\sigma_1, \infty}(\widetilde{\Omega}_T)$ implies the same regularity of the momentum: $\rho^\varepsilon\mathbf u^\varepsilon\in  B_3^{\sigma_1, \infty}(\widetilde{\Omega}_T)$.
psilon
\item
As pointed out in \cite{chen2024}, the condition $\gamma\geq2$ could be relaxed in view of~\cite{ak}. 

\item
The regularity in the interior is in line with the positive half of Onsager's conjecture, established by Constantin--E--Titi~\cite{c}: the H\"{o}lder $C^{0,1/3+}$-continuity of the velocity ensures the energy conservation for weak solutions to the incompressible Euler equations. However, in general, the absence of anomalous dissipation does not necessarily guarantee the energy conservation in the vanishing viscosity limit.
\end{enumerate}

The remaining parts of the paper are structured as follows. Section~\ref{sec: prelim} presents the mathematical preliminaries; in particular, the commutator estimates and the boundary layer foliation. Section~\ref{sec: proof of main theorem} is devoted to the proof of the Kato-type criterion for MHD, \textit{i.e.}, the first half of Theorem~\ref{thm: main}. Finally, in Section~\ref{sec: proof, convergence} we prove the second half of Theorem~\ref{thm: main}, namely the convergence of solutions in the inviscid limit.

\section{Preliminaries}\label{sec: prelim}
\subsection{Commutator estimate and Hardy's inequality}
\hspace{1.5em}

For $\Omega \subset \R$, recall the Besov space \(B^{\alpha,\infty}_p(\Omega)\) for \(p \in [1, \infty]\) and \(\alpha \in (0, 1)\), which consists of functions $f:\Omega\to\R$ with finite Besov norm:
\begin{equation}\label{2.2}
\|f\|_{B^{\alpha,\infty}_p(\Omega)} := \|f\|_{L^p(\Omega)} + [f]_{B^{\alpha,\infty}_p(\Omega)},
\end{equation}
where the H\"older-type semi-norm is defined as
\begin{equation}\label{2.3}
[f]_{B^{\alpha,\infty}_p(\Omega)} := \sup_{|y|>0} \frac{\|f(\cdot + y) - f(\cdot)\|_{L^p(\Omega \cap (\Omega - \{y\}))}}{|y|^\alpha}.
\end{equation}

Let $\eta_\e$ be a standard mollifier of width $\e$, namely that $\eta_\e(x) = \e^{-3}\eta(x/\e)$ with a smooth, even kernel $\eta$ compactly supported in the unit ball with $\int_{\R}\eta\,\dd x=1$. Throughout this paper, denote by $\overline{f}_\e$ the mollification of $f: \Omega \to\mathbb{R}$ with $\eta_\e$: 
\begin{equation}\label{2.1}
\overline{f}_\varepsilon(x) := \int_{B_\varepsilon(0)} f(x - y)\eta_\varepsilon(y)\,\dd y, \qquad \text{for any } x \in \Omega^\varepsilon,
\end{equation}
where $$\Omega^\e := \{x \in \Omega: {\rm dist}(x,\p\Omega)>\e\}.$$ Then, for any \(f \in B_r^{\alpha,\infty}(\Omega)\) with \(r \in [1,\infty]\), it holds that
\begin{equation}\label{2.4}
\left\|\nabla \overline{f}_\varepsilon\right\|_{L^r(\Omega^\varepsilon)} \leq C \varepsilon^{\alpha-1} \|f\|_{B_r^{\alpha,\infty}(\Omega)}
\end{equation}
and
\begin{equation}\label{2.5}
\left\|\overline{f}_\varepsilon - f\right\|_{L^r(\Omega^\varepsilon)} \leq C\varepsilon^{\alpha} \|f\|_{B_r^{\alpha,\infty}(\Omega)}.
\end{equation}

In this paper, we shall frequently use the following commutator estimate and Hardy's inequality: \begin{Lemma}[\cite{f}] \label{lem:commutator}
Given \( f \in B_{r_1}^{\alpha,\infty}(\Omega) \), \( g \in B_{r_2}^{\alpha,\infty}(\Omega) \) with \( 1 \leq r, r_1, r_2 < \infty \) and $\frac{1}{r_1} + \frac{1}{r_2} = \frac{1}{r}$. Then, for some constant $C$ independent of $\e$, we have that
\begin{equation}\label{2.7}
\left\| \overline{(f \otimes g)}_{\varepsilon} - \bar{f}_{\varepsilon} \otimes \bar{g}_{\varepsilon} \right\|_{L^r(\Omega^{\varepsilon})} \leq C \varepsilon^{2\alpha} \| f \|_{B_{r_1}^{\alpha,\infty}(\Omega)} \| g \|_{B_{r_2}^{\alpha,\infty}(\Omega)}.
\end{equation}
\end{Lemma}

\begin{Lemma}[\cite{kufner77}]\label{lem: hardy}
Let \(p \in [1, \infty)\) and \(f \in W_0^{1,p}(\Omega)\) such that the following inequality holds:
\begin{equation}\label{2.8}
\left\| \frac{f}{\operatorname{dist}(\cdot, \partial\Omega)} \right\|_{L^p(\Omega)} \leq C \|\nabla f\|_{L^p(\Omega)},
\end{equation}
where the constant \(C\) depends on \(p\) and \(\Omega\).
\end{Lemma}

\subsection{Boundary layer foliation}
\label{subsec:viscous_estimation}

In order to obtain delicate quantitative estimates in the boundary layer, for subsequence developments we make crucial use of the domain foliation strategy as in Chen--Liang--Wang~\cite{chen2022}.

Given regularity parameters $\sigma_1, \sigma_2 \in \left(\frac{1}{3}, \frac{1}{2}\right]$ as in \eqref{1.5}; denote $\sigma = \min\{\sigma_1, \sigma_2\}$. We define a sequence $\{\beta_n^*\}_{n\ge 0}$ as follows:
\begin{equation}\label{2.1}
\begin{cases}
\beta_0^* = 0, \\
 \beta_n^* = \frac{1}{2(1-\sigma)} \left(1 + \frac{1}{3}\beta_{n-1}^*\right), \quad n = 1, 2, 3, \cdots.
\end{cases} 
\end{equation}
It is clear that
\begin{equation}\label{2.2}
\beta_\infty^* := \lim_{n \to \infty} \beta_n^* = \frac{3}{5-6\sigma} > 1,
\end{equation}
so there exists a finite integer $N = N(\sigma) \ge 1$ such that
\begin{equation}\label{2.3}
0 = \beta_0^* < \beta_1^* < \beta_2^* < \cdots < \beta_{N-1}^* \le 1 < \beta_N^*.
\end{equation}

By perturbing $\beta_n^*$ slightly, one can choose a nearby sequence $\{\beta_n\}_{n=0}^N$:
\begin{equation}\label{2.12}
0 = \beta_0 < \beta_n < \frac{1}{2(1-\sigma)} \left(1 + \frac{1}{3}\beta_{n-1}\right) \quad \text{for } n \ge 1,
\end{equation}
and
\begin{equation}\label{2.13}
\beta_0 < \beta_1 (=\beta) < \beta_2 < \cdots < \beta_{N-1} < 1 = \beta_N.
\end{equation}

Based on \eqref{2.4}--\eqref{2.5} and (1.14), the physical domain $\Omega$ is partitioned into subregions as follows:
\begin{equation}\label{2.14}
\begin{cases}
    V_1 = \Omega^{\varepsilon^{\beta_1}}, \\
    V_n = \Omega^{\varepsilon^{\beta_n}} \setminus \Omega^{\varepsilon^{\beta_{n-1}} + \varepsilon^{\beta_n}} \quad \quad (2 \le n \le N)\\
    V_{N+1} = \Omega \setminus \left(\bigcup_{n=1}^N V_n\right) =: \Gamma_\varepsilon. 
\end{cases}
\end{equation}
For $\varepsilon > 0$ sufficiently small, one has
\begin{equation}\label{2.15}
V_k \cap V_m = \emptyset \qquad \text{whenever } |k - m| \ge 2.
\end{equation}
Furthermore, for some constant $C$ depending only on the geometry of $\Omega$,
\begin{equation}\label{2.16}
\begin{cases}
\operatorname{meas} V_n \le C \varepsilon^{\beta_{n-1}}, \\ \operatorname{meas} (V_n \cap V_{n+1}) \le C \varepsilon^{\beta_{n+1}},\\ \operatorname{meas} (V_n \cap V_{n+1}^c) \le C \varepsilon^{\beta_{n-1}}.
\end{cases}
\end{equation}
where $V^c$ denotes the complement of $V$ and ${meas}$ denotes the Lebesgue measure. It follows that
\begin{equation}\label{2.18}
\operatorname{meas} (V_k \cap V_m) \le
\begin{cases}
C\delta^{\min\{\beta_k,\beta_m\}} & \text{if } |k-m| = 1, \\
0 & \text{if } |k-m| > 1,
\end{cases}
\end{equation}
whenever $\e>0$ is small enough.  

Let  $\{\xi_n\}_{n=1}^{N}$ be a smooth partition of unity subordinate to $\{V_n\}_{n=1}^{N}$. Then
\begin{equation}\label{2.19}
\operatorname{spt}\xi_n \subset V_n, \qquad 0 \le \xi_n \le 1, \qquad  \sum_{n=1}^{N} \xi_n(x) = 1 \text{ for all } x \in \Omega.
\end{equation}

Observe by~\eqref{2.19} and \eqref{2.15} that 
\begin{equation}\label{2.20}
\nabla(\xi_n + \xi_{n+1})^2 = 0 \qquad \text{on } V_n \cap V_{n+1},
\end{equation}
namely that the gradient of $\xi$'s vanishes identically on overlapping regions. We shall repeatedly use this simple identity.

\begin{figure}[htbp]
  \centering
  \includegraphics[width=1\textwidth]
  {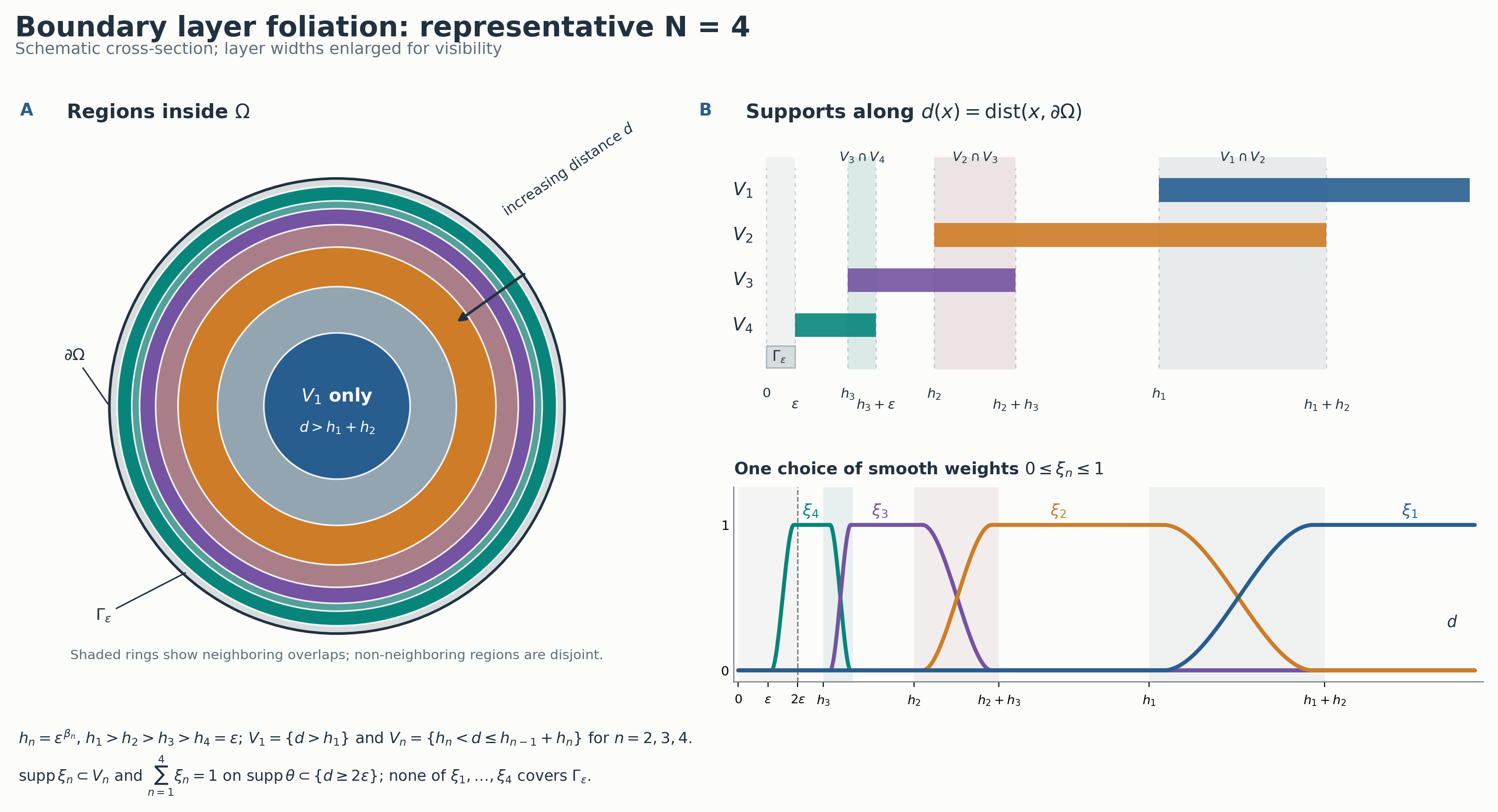}
  \caption{A schematic picture explaining the construction of $\{V_n\}$ and $\{\xi_n\}$}
\end{figure}


\section{Absence of anomalous dissipation}\label{sec: proof of main theorem}

This section is devoted to the proof of the first half of Theorem~\ref{thm: main}, namely the absence of anomalous dissipation under the Kato-type condition. That is, we show that~\eqref{1.8} implies~\eqref{1.7}

In the proof, we drop the superscript $\varepsilon$ and denote $(\rho^\varepsilon, \mathbf{u}^\varepsilon, \mathbf{B}^\varepsilon)$ by $(\rho, \mathbf{u}, \mathbf{B})$ for convenience. Throughout the proof, $C > 0$ denotes a generic positive constant that may depend on $T, \widehat{\rho}, \gamma, \lambda$, and the domain $\Omega$, but is strictly independent of $\varepsilon$.

\subsection{An identity for the approximate energy}
We first derive an identity~\eqref{3.20} for the energy $\mathcal{E}^\varepsilon(t)
=\int_\Omega
\left(\frac{1}{2}\rho^\varepsilon |\mathbf{u}^\varepsilon|^2
+\frac{(\rho^\varepsilon)^\gamma}{\gamma-1}
+\frac{1}{2}|\mathbf{B}^\varepsilon|^2\right)\,\dd x$ with suitable spacetime cutoffs and regularisations. It is the starting point for all subsequent analyses.

For ease of future reference, we denote for $n \in \mathbb{N}$, $T>0$, and $\e>0$ sufficiently small that
\begin{align*}
    &V_T^n := (V_n \cap V_{n+1}^c) \times (\varepsilon^{\beta_1}, T - \varepsilon^{\beta_1}),\\
    &V_T^{n,n+1} = (V_n \cap V_{n+1}) \times (\varepsilon^{\beta_1}, T - \varepsilon^{\beta_1}).
\end{align*}
Here $\beta_1=\beta$ is as in \S\ref{subsec:viscous_estimation} above. 

Recall the definition of $V_n$ in \eqref{2.14}. The following mollifications for a given $f: \Omega \to \mathbb{R}$ are well-defined: for each $n \in \{1, 2, \dots, N\}$, 
\begin{equation}\label{3.1}
    \overline{f}_n(x,t) := \overline{f}_{\varepsilon^{\beta_n}}(x,t) \equiv
    \begin{cases}
        \int_0^T \int_\Omega \eta_{\varepsilon^{\beta_n}}(x-y, t-s)f(y,s)\,\dd y\,\dd s \qquad\text{for } (x,t) \in V_T^n,\\
        \int_0^T \int_\Omega \eta_{\varepsilon^{\beta_{n+1}}}(x-y, t-s)f(y,s)\,\dd y\,\dd s \qquad\text{for } (x,t) \in V_T^{n,n+1}.
    \end{cases}
\end{equation}
In particular, the mollification width varies as $x$ reaches different annuli $V_n$. Note that
\begin{equation}\label{3.2}
    \overline{f}_n = \overline{f}_{n+1} \quad \text{on } V_T^{n,n+1}.
\end{equation}

Denote by $\partial$ either $\partial_{x_i}$ or $\partial_t$. For all $p \in [1, \infty)$ one has that
\begin{equation}\label{3.3}
    \|\partial \overline{f}_n\|_{L^p} \leqslant
    \begin{cases}
        C\varepsilon^{\beta_n(\alpha-1)}\|f\|_{B_p^{\alpha, \infty}(\Omega_T)}, & (x,t) \in V_T^n, \\[2ex]
        C\varepsilon^{\beta_{n+1}(\alpha-1)}\|f\|_{B_p^{\alpha, \infty}(\Omega_T)}, & (x,t) \in V_T^{n,n+1},
    \end{cases}
\end{equation}
as well as
\begin{equation}\label{3.4}
    \|\overline{f}_n - f\|_{L^p} \leqslant
    \begin{cases}
        C\varepsilon^{\beta_n\alpha}\|f\|_{B_p^{\alpha, \infty}(\Omega_T)}, & (x,t) \in V_T^n, \\[2ex]
        C\varepsilon^{\beta_{n+1}\alpha}\|f\|_{B_p^{\alpha, \infty}(\Omega_T)}, & (x,t) \in V_T^{n,n+1}.
    \end{cases}
\end{equation}

Applying the above mollification term-wise to the MHD system~\eqref{1.1} leads to the regularised equations, which holds \textit{a.e.} in $V_n \times (\varepsilon^{\beta_1}, T - \varepsilon^{\beta_1})$ for each $n \in \{1, 2, \dots, N\}$:
\begin{equation}\label{3.5}
\begin{aligned}
    \partial_t \overline{(\rho \mathbf{u})}_n + \operatorname{div} \overline{(\rho \mathbf{u} \otimes \mathbf{u})}_n + \nabla \overline{(\rho^\gamma)}_n &= \varepsilon \overline{\Delta \mathbf{u}}_n + \overline{(\operatorname{curl} \mathbf{B} \times \mathbf{B})}_n, \\[1ex]
    \partial_t \overline{\mathbf{B}}_n - \operatorname{curl} \overline{(\mathbf{u} \times \mathbf{B})}_n &= - \varepsilon \operatorname{curl} (\operatorname{curl} \overline{\mathbf{B}}_n).
\end{aligned}
\end{equation}
We multiply \eqref{3.5} by $\xi_n$, the partition of unity in \eqref{2.19}, and sum over $n$ to deduce that  
\begin{equation}\label{eq:3.7}
\begin{aligned}
    \sum_{n=1}^N \xi_n \bigl[ & \partial_t (\overline{\rho}_n\,\overline{\mathbf{u}}_n)
     + \operatorname{div}(\overline{(\rho\mathbf{u})}_n \otimes \overline{\mathbf{u}_n})
     + \nabla (\bar{\rho}_n)^\gamma \bigr]\\
    = \sum_{n=1}^N \xi_n \bigl[ &\partial_t \bigl(
        \overline{\rho}_n\,\overline{\mathbf{u}}_n - \overline{(\rho \mathbf{u})}_n \bigr)
    + \operatorname{div}(\overline{(\rho\mathbf{u})}_n \otimes \overline{\mathbf{u}_n}
        - \overline{(\rho \mathbf{u} \otimes \mathbf{u})}_n \bigr)\\
    &+ \nabla \bigl( (\overline{\rho}_n)^\gamma - \overline{(\rho^\gamma)}_n \bigr)
     + \varepsilon \Delta \overline{\mathbf{u}}_n
     + (\overline{\operatorname{curl} \mathbf{B} \times \mathbf{B}})_n \bigr]
\end{aligned}
\end{equation}
and
\begin{equation}\label{eq:3.8}
\begin{aligned}
    \sum_{n=1}^N \xi_n \bigl[ \partial_t \overline{\mathbf{B}}_n
    - \operatorname{curl}(\overline{\mathbf{u} \times \mathbf{B}})_n \bigr]
    = -\, \varepsilon \sum_{n=1}^N \xi_n
        \operatorname{curl}\bigl(\operatorname{curl}\overline{\mathbf{B}}_n\bigr).
\end{aligned}
\end{equation}

We further apply a global spacetime cutoff: choose $\psi(t) \in C_0^1(0, T)$ and $\theta \in C^\infty_0(\Omega)$ such that
\begin{equation}\label{3.9}
    0 \leqslant \theta \leqslant 1 \text{ in } \Omega, \quad \theta = 1 \text{ on } \Omega^{4\varepsilon}, \quad \theta = 0 \text{ on } \Gamma_{2\varepsilon}, \quad \text{and }|\nabla \theta| \leqslant \frac{4}{\varepsilon}.
\end{equation}
Recall that $\nabla\theta = 0$ for $x \notin \Omega^{2\varepsilon} \cap \Gamma_{4\varepsilon}$. Then, by the definition of $V_N$, it is evident that $V_N$ contains the entire region where $\nabla\theta \neq 0$; equivalently, $\operatorname{supp}(\nabla\theta) \subset V_N$. As a result, all the integrals involving $\nabla\theta$ can be restricted to $V_N$.

Then, from \eqref{eq:3.7} and \eqref{eq:3.8} we deduce that
\begin{align}\label{3.10}
    &\int_0^T \psi(t) \int_\Omega \theta(x) \bigg(\sum_{m=1}^N \xi_m \overline{\mathbf{u}}_m\bigg) \cdot \bigg( \sum_{n=1}^N \xi_n \partial_t (\overline{\rho}_n\,\overline{\mathbf{u}}_n) + \sum_{n=1}^N \xi_n \operatorname{div}(\overline{(\rho\mathbf{u})}_n \otimes \overline{\mathbf{u}}_n) \bigg) \,\dd x\,\dd t \nonumber\\
    &+ \int_0^T \psi(t) \int_\Omega \theta(x) \bigg(\sum_{m=1}^N \xi_m \overline{\mathbf{u}}_m\bigg) \cdot \bigg( \sum_{n=1}^N \xi_n \nabla\big((\overline{\rho}_n)^\gamma\big) \bigg) \,\dd x\,\dd t \nonumber\\
    &+ \int_0^T \psi(t) \int_\Omega \theta(x) \bigg(\sum_{m=1}^N \xi_m \overline{\mathbf{B}}_m\bigg) \cdot \bigg(\sum_{n=1}^N \xi_n \partial_t \overline{\mathbf{B}}_n\bigg) \,\dd x\,\dd t \nonumber\\
    &= r_1^\varepsilon + r_2^\varepsilon + r_3^\varepsilon + \varepsilon \int_0^T \psi(t) \int_\Omega \theta(x) \bigg(\sum_{m=1}^N \xi_m \overline{\mathbf{u}}_m\bigg) \cdot \bigg(\sum_{n=1}^N \xi_n \Delta \overline{\mathbf{u}}_n\bigg) \,\dd x\,\dd t \nonumber\\
    &\quad + \int_0^T \psi(t) \int_\Omega \theta(x) \bigg(\sum_{m=1}^N \xi_m \overline{\mathbf{u}}_m\bigg) \cdot \bigg(\sum_{n=1}^N \xi_n (\overline{\operatorname{curl} \mathbf{B} \times \mathbf{B}})_n\bigg) \,\dd x\,\dd t \nonumber\\
    &\quad + \int_0^T \psi(t) \int_\Omega \theta(x) \bigg(\sum_{m=1}^N \xi_m \overline{\mathbf{B}}_m\bigg) \cdot \bigg(\sum_{n=1}^N \xi_n \operatorname{curl}(\overline{(\mathbf{u} \times \mathbf{B})}_n)\bigg) \,\dd x\,\dd t \nonumber\\
    &\quad - \varepsilon \int_0^T \psi(t) \int_\Omega \theta(x) \bigg(\sum_{m=1}^N \xi_m \overline{\mathbf{B}}_m\bigg) \cdot \bigg(\sum_{n=1}^N \xi_n \operatorname{curl}(\operatorname{curl}\overline{\mathbf{B}}_n)\bigg) \,\dd x\,\dd t,
\end{align}
with the commutator terms:
\begin{align*}
    r_1^\varepsilon &= \int_0^T \psi(t) \int_\Omega \theta(x) \bigg(\sum_{m=1}^N \xi_m \overline{\mathbf{u}}_m\bigg) \cdot \bigg( \sum_{n=1}^N \xi_n \partial_t \big( (\overline{\rho}_n\,\overline{\mathbf{u}}_n) - \overline{(\rho\mathbf{u})}_n \big) \bigg) \,\dd x\,\dd t, \\
    r_2^\varepsilon &= \int_0^T \psi(t) \int_\Omega \theta(x) \bigg(\sum_{m=1}^N \xi_m \overline{\mathbf{u}}_m\bigg) \cdot \bigg( \sum_{n=1}^N \xi_n \operatorname{div} \big( (\overline{(\rho\mathbf{u})}_n \otimes \overline{\mathbf{u}}_n) - \overline{(\rho\mathbf{u} \otimes \mathbf{u})}_n \big) \bigg) \,\dd x\,\dd t, \\
    r_3^\varepsilon &= \int_0^T \psi(t) \int_\Omega \theta(x) \bigg(\sum_{m=1}^N \xi_m \overline{\mathbf{u}}_m\bigg) \cdot \bigg( \sum_{n=1}^N \xi_n \nabla \big( (\overline{\rho}_n)^\gamma - \overline{(\rho^\gamma)}_n \big) \bigg) \,\dd x\,\dd t.
\end{align*}

In the sequel, we shall rewrite~\eqref{3.10} by expanding the terms on its left-hand side.

\begin{itemize}
    \item 
For the first term on the left-hand side of~\eqref{3.10}, we note by $\xi_n\xi_m \equiv 0$ for $|n-m|\geq 2$ that 
\begin{equation}\label{eq:3.12}
\begin{aligned}
    &\int_0^T \psi \int_\Omega \theta \bigg(\sum_{m=1}^N \xi_m \overline{\mathbf{u}}_m\bigg) \cdot \partial_t \bigg(\sum_{n=1}^N \xi_n (\overline{\rho}_n\,\overline{\mathbf{u}}_n)\bigg) \,\dd x\,\dd t \\
    &= \sum_{|m-n|\le 1} \int_0^T \psi \int_\Omega \theta \xi_m \xi_n \overline{\mathbf{u}}_m \cdot \partial_t (\overline{\rho}_n\,\overline{\mathbf{u}}_n) \,\dd x\,\dd t \\
    &= \sum_{|m-n|\le 1} \int_0^T \psi \int_\Omega \theta \xi_m \xi_n \left[ \partial_t \left( \overline{\rho}_{\tilde{n}} \frac{|\overline{\mathbf{u}}_{\tilde{n}}|^2}{2} \right) + \partial_t \overline{\rho}_{\tilde{n}} \frac{|\overline{\mathbf{u}}_{\tilde{n}}|^2}{2} \right] \,\dd x\,\dd t.
\end{aligned}
\end{equation}
Here and hereafter, we introduce the notation as in \cite{chen2022}:
\begin{equation}\label{eq:3.13}
    \tilde{n} =
    \begin{cases}
        \max\{m, n\} &\qquad \text{if } |m - n| = 1, \\
        n &\qquad \text{if } |m - n| = 0.
    \end{cases}
\end{equation}

To see the final equality in \eqref{eq:3.12}, note that on the support of $\xi_m\xi_n$ with $|m-n|\le 1$, we have $\overline{\mathbf{u}}_m \cdot \partial_t (\overline{\rho}_n\,\overline{\mathbf{u}}_n) = \overline{\mathbf{u}}_{\tilde{n}} \cdot \partial_t (\overline{\rho}_{\tilde{n}}\,\overline{\mathbf{u}}_{\tilde{n}})$ in light of~\eqref{3.2}. We further compute
\begin{align*}
\overline{\mathbf{u}}_{\tilde{n}} \cdot \partial_t (\overline{\rho}_{\tilde{n}}\,\overline{\mathbf{u}}_{\tilde{n}}) &= \overline{\mathbf{u}}_{\tilde{n}} \cdot \left( (\partial_t \overline{\rho}_{\tilde{n}}) \overline{\mathbf{u}}_{\tilde{n}} + \overline{\rho}_{\tilde{n}} \partial_t \overline{\mathbf{u}}_{\tilde{n}} \right) \\
    &= \left[ (\partial_t \overline{\rho}_{\tilde{n}}) \frac{|\overline{\mathbf{u}}_{\tilde{n}}|^2}{2} + \overline{\rho}_{\tilde{n}} \partial_t \left(\frac{|\overline{\mathbf{u}}_{\tilde{n}}|^2}{2}\right) \right] + (\partial_t \overline{\rho}_{\tilde{n}}) \frac{|\overline{\mathbf{u}}_{\tilde{n}}|^2}{2} \\
    &= \partial_t \left( \overline{\rho}_{\tilde{n}} \frac{|\overline{\mathbf{u}}_{\tilde{n}}|^2}{2} \right) + \partial_t \overline{\rho}_{\tilde{n}} \frac{|\overline{\mathbf{u}}_{\tilde{n}}|^2}{2}.
\end{align*}

By a similar argument, we obtain
\begin{equation}\label{eq:3.14}
\begin{aligned}
    &\int_0^T \psi \int_\Omega \theta \bigg(\sum_{m=1}^N \xi_m \overline{\mathbf{u}}_m\bigg) \cdot \bigg( \sum_{n=1}^N \xi_n \operatorname{div}\big((\overline{\rho \mathbf{u}})_n \otimes \overline{\mathbf{u}}_n\big) \bigg)\,\dd x\,\dd t \\
    &\qquad = \sum_{|m-n|\le 1} \int_0^T \psi \int_\Omega \theta \xi_m \xi_n \left[ \operatorname{div}\left( (\overline{\rho \mathbf{u}})_{\tilde{n}} \frac{|\overline{\mathbf{u}}_{\tilde{n}}|^2}{2} \right) + \operatorname{div}(\overline{\rho \mathbf{u}})_{\tilde{n}} \frac{|\overline{\mathbf{u}}_{\tilde{n}}|^2}{2} \right]\,\dd x\,\dd t.
\end{aligned}
\end{equation}

To proceed, in view of \eqref{3.1}, \eqref{3.2}, and \eqref{eq:3.13}, we mollify the mass equation in the MHD system~\eqref{1.2} to deduce for each $n$ that
\begin{equation}\label{mollification of mass eq}
    0 = \partial_t \bar{\rho}_n + \operatorname{div}(\overline{\rho \mathbf{u}})_n \quad \text{a.e. in } V_n \times (0, T).
\end{equation}
This together with \eqref{eq:3.12} and \eqref{eq:3.14} yields that
\begin{align}
    &\int_0^T \psi \int_\Omega \theta \bigg(\sum_{m=1}^N \xi_m \overline{\mathbf{u}}_m\bigg) \cdot \partial_t \bigg(\sum_{n=1}^N \xi_n (\overline{\rho}_n\,\overline{\mathbf{u}}_n)\bigg)\,\dd x\,\dd t \nonumber \\
    &\quad + \int_0^T \psi \int_\Omega \theta \bigg(\sum_{m=1}^N \xi_m \overline{\mathbf{u}}_m\bigg) \cdot \sum_{n=1}^N \xi_n \operatorname{div}\big((\overline{\rho \mathbf{u}})_n \otimes \overline{\mathbf{u}}_n\big) \,\dd x\,\dd t\nonumber \\
    &= \sum_{|m-n|\le 1} \int_0^T \psi \int_\Omega \theta \xi_m \xi_n \left[ \partial_t \left( \overline{\rho}_{\tilde{n}} \frac{|\overline{\mathbf{u}}_{\tilde{n}}|^2}{2} \right) + \operatorname{div}\left( (\overline{\rho \mathbf{u}})_{\tilde{n}} \frac{|\overline{\mathbf{u}}_{\tilde{n}}|^2}{2} \right) \right]\,\dd x\,\dd t \nonumber \\
    &= - \sum_{|m-n|\le 1} \int_0^T \psi' \int_\Omega \theta \xi_m \xi_n \overline{\rho}_{\tilde{n}} \frac{|\overline{\mathbf{u}}_{\tilde{n}}|^2}{2}\,\dd x\,\dd t \nonumber \\
    &\quad \underbrace{- \sum_{|m-n|\le 1} \int_0^T \psi \int_\Omega \nabla \theta \xi_m \xi_n \cdot (\overline{\rho \mathbf{u}})_{\tilde{n}} \frac{|\overline{\mathbf{u}}_{\tilde{n}}|^2}{2}\,\dd x\,\dd t}_{=:\,J_1} \underbrace{- \sum_{|m-n|\le 1} \int_0^T \psi \int_\Omega \theta \nabla(\xi_m \xi_n) \cdot (\overline{\rho \mathbf{u}})_{\tilde{n}} \frac{|\overline{\mathbf{u}}_{\tilde{n}}|^2}{2}\,\dd x\,\dd t}_{=:\,J_2}. \label{3.15}
\end{align}

\item 
For the second term on the left-hand side of \eqref{3.10},  similar arguments lead to
\begin{align}\label{eq:3.16}
&\int_0^T \psi \int_\Omega \theta \bigg(\sum_{m=1}^N \xi_m \overline{\mathbf{u}}_m\bigg) \cdot \bigg(\sum_{n=1}^N \xi_n \nabla\big((\bar{\rho}_n)^\gamma\big)\bigg)\,\dd x\,\dd t\nonumber\\
&\qquad = \sum_{|m-n|\le 1} \int_0^T \psi \int_\Omega \theta \xi_m \xi_n \overline{\mathbf{u}}_m \cdot \nabla\big((\bar{\rho}_n)^\gamma\big)\,\dd x\,\dd t.
\end{align}
Observe that the mollified mass equation~\eqref{mollification of mass eq} can be expressed as 
\[
    \partial_t \bar{\rho}_n + \operatorname{div}(\bar{\rho}_n \overline{\mathbf{u}}_n)
    = \operatorname{div}\big(\bar{\rho}_n \overline{\mathbf{u}}_n - (\overline{\rho \mathbf{u}})_n\big),
\]
and hence 
\begin{align}\label{3.17}
    \partial_t (\bar{\rho}_n)^\gamma + (\gamma-1)(\bar{\rho}_n)^\gamma \operatorname{div} \overline{\mathbf{u}}_n + \operatorname{div}\big((\bar{\rho}_n)^\gamma \overline{\mathbf{u}}_n\big) &= \gamma(\bar{\rho}_n)^{\gamma-1} \operatorname{div}\big(\bar{\rho}_n \overline{\mathbf{u}}_n - (\overline{\rho \mathbf{u}})_n\big)
\end{align}
by multiplication with $\gamma(\bar{\rho}_n)^{\gamma-1}$ on both sides. Applying~\eqref{3.17} and integration by parts, the right-hand side of \eqref{eq:3.16} thus becomes
\begin{align}
    &\sum_{|m-n|\le 1} \int_0^T \psi \int_\Omega \theta \xi_m \xi_n \overline{\mathbf{u}}_n \cdot \nabla \big((\bar{\rho}_n)^\gamma\big) \,\dd x\,\dd t \nonumber\\
    &\qquad = J_3+J_4+J_5 - \frac{1}{\gamma-1} \sum_{|m-n|\le 1} \int_0^T \psi' \int_\Omega \theta \xi_m \xi_n (\bar{\rho}_n)^\gamma \,\dd x\,\dd t,
    \label{eq:3.18}
\end{align}
where 
\begin{align*}
    &J_3 := - \sum_{|m-n|\le 1} \int_0^T \psi \int_\Omega \nabla(\theta \xi_m \xi_n) \cdot \overline{\mathbf{u}}_n (\bar{\rho}_n)^\gamma \,\dd x\,\dd t,\\
    &J_4:= \frac{1}{\gamma-1} \sum_{|m-n|\le 1} \int_0^T \psi \int_\Omega \theta \xi_m \xi_n \operatorname{div}\big((\bar{\rho}_n)^\gamma \overline{\mathbf{u}}_n\big) \,\dd x\,\dd t,\\
    &J_5:=- \frac{\gamma}{\gamma-1} \sum_{|m-n|\le 1} \int_0^T \psi \int_\Omega \theta \xi_m \xi_n (\bar{\rho}_n)^{\gamma-1} \operatorname{div}\big(\bar{\rho}_n \overline{\mathbf{u}}_n - (\overline{\rho \mathbf{u}})_n\big)\,\dd x\,\dd t.
\end{align*}

\item 
Next, for the final term on the left-hand side of \eqref{3.10}, integration by parts in time yields
\begin{align}
    &\int_0^T \psi \int_\Omega \theta \left(\sum_{n=1}^N \xi_n \overline{\mathbf B}_n\right) \cdot \left(\sum_{n=1}^N \xi_n \partial_t \overline{\mathbf B}_n\right) \,\dd x\,\dd t \nonumber\\
    &\qquad = \sum_{|m-n|\leqslant 1} \int_0^T \psi \int_\Omega \xi_m \xi_n \, \partial_t\left(\frac{1}{2}|\overline{\mathbf B}_n|^2\right) \,\dd x\,\dd t \nonumber \\
    &\qquad= -\sum_{|m-n|\leqslant 1} \int_0^T \psi' \int_\Omega \xi_m \xi_n \left(\frac{1}{2}|\overline{\mathbf B}_n|^2\right) \dd x\,\dd t.
    \label{3.19}
\end{align}
\end{itemize}

Therefore, summarising \eqref{3.15}--\eqref{3.19}, we may recast \eqref{3.10} into the following:
\begin{align}\label{3.20}
    &- \sum_{|m-n|\le 1} \int_0^T \psi' \int_\Omega \theta \xi_m \xi_n \bigg( \overline{\rho}_{\tilde{n}} \frac{|\overline{\mathbf{u}}_{\tilde{n}}|^2}{2} + \frac{(\overline{\rho}_{\tilde{n}})^\gamma}{\gamma-1} + \frac{|\overline{\mathbf{B}}_{\tilde{n}}|^2}{2} \bigg) \,\dd x\,\dd t \nonumber\\
    &\qquad= r_1^\varepsilon + r_2^\varepsilon + r_3^\varepsilon  - \sum_{j=1}^5 J_j \nonumber\\
    &\qquad\quad + \varepsilon \int_0^T \psi(t) \int_\Omega \theta(x) \bigg(\sum_{m=1}^N \xi_m \overline{\mathbf{u}}_m\bigg) \cdot \bigg(\sum_{n=1}^N \xi_n \Delta \overline{\mathbf{u}}_n\bigg) \,\dd x\,\dd t \nonumber\\
    &\qquad\quad + \int_0^T \psi(t) \int_\Omega \theta(x) \bigg(\sum_{m=1}^N \xi_m \overline{\mathbf{B}}_m\bigg) \cdot \bigg(\sum_{n=1}^N \xi_n \operatorname{curl}(\overline{\mathbf{u} \times \mathbf{B}})_n\bigg) \,\dd x\,\dd t \nonumber\\
    &\qquad\quad - \varepsilon \int_0^T \psi(t) \int_\Omega \theta(x) \bigg(\sum_{m=1}^N \xi_m \overline{\mathbf{B}}_m\bigg) \cdot \bigg(\sum_{n=1}^N \xi_n \operatorname{curl}(\operatorname{curl}\overline{\mathbf{B}}_n)\bigg) \,\dd x\,\dd t \nonumber\\
    &\qquad\quad + \int_0^T \psi(t) \int_\Omega \theta(x) \bigg(\sum_{m=1}^N \xi_m \overline{\mathbf{u}}_m\bigg) \cdot \bigg(\sum_{n=1}^N \xi_n (\overline{\operatorname{curl}\mathbf{B} \times \mathbf{B}})_n\bigg) \,\dd x\,\dd t.
\end{align}
See below~\eqref{3.10} for the expressions of $r_1^\e$, $r_2^\e$, and $r_3^\e$.

\subsection{Proof of the Kato-type criterion} 
In this subsection, we prove~\eqref{1.7}, reproduced below:
\begin{align*}
     \lim_{\varepsilon \to 0} \varepsilon \int_0^T \int_\Omega \left( |\nabla \mathbf{u}^\varepsilon|^2 + |\nabla \mathbf{B}^\varepsilon|^2 \right) \,\dd x \,\dd t = 0,
\end{align*}
provided that \eqref{1.8} holds, namely the same result with $\int_\Omega$ replaced by $\int_{\Gamma_{8\e}}$. The key is to exploit the left-most term in \eqref{3.20} above.

First of all, following almost verbatim the arguments in~\cite{chen2024}, we arrive at the following. 

\begin{proposition}\label{prop: 3.1}
    Under the assumptions of Theorem~\ref{thm: main}, it holds that
\begin{equation}\label{3.21}
    \lim_{\varepsilon \to 0} \left| r_1^\varepsilon + r_2^\varepsilon + r_3^\varepsilon
    - \sum_{j=1}^5 J_j
    + \varepsilon \int_0^T \psi(t) \int_\Omega \theta(x)
    \bigg(\sum_{m=1}^N \xi_m \overline{\mathbf{u}}_m\bigg)
    \cdot \bigg(\sum_{n=1}^N \xi_n \Delta \overline{\mathbf{u}}_n\bigg)
    \,\dd x\,\dd t \right| = 0.
\end{equation}
\end{proposition}

\vspace{1em}

Now, let us turn to the magnetic diffusion term.
\begin{proposition}\label{prop: 3.2}
Under the assumptions of Theorem~\ref{thm: main}, it holds that
\begin{equation}\label{3.22}
    \lim_{\varepsilon \to 0} \left| \varepsilon \int_0^T \psi(t) \int_\Omega \theta(x) \left( \sum_{m=1}^N \xi_m \overline{\boldsymbol{B}_m} \right) \cdot \left( \sum_{n=1}^N \xi_n \operatorname{curl} (\operatorname{curl} \overline{\boldsymbol{B}_n}) \right) \,\dd x\,\dd t \right| = 0.
\end{equation}
\end{proposition}

\begin{proof}
Thanks to the choice of the partition of unity $\{\xi_n\}$ in~\eqref{2.19}, the nonzero contribution to the sum over $m,n$ arises  from $\{(m,n)\in \{1,\ldots N\}^2:\,|m-n|\leq 1\}$ only. Then, in view of the vector identity $\operatorname{div}(\boldsymbol{A} \times \boldsymbol{C}) = \boldsymbol{C} \cdot \operatorname{curl} \boldsymbol{A} - \boldsymbol{A} \cdot \operatorname{curl} \boldsymbol{C}$ and integration by parts, we deduce that
\begin{align}
    & \varepsilon \int_0^T \psi(t) \int_\Omega \theta(x) \left( \sum_{m=1}^N \xi_m \overline{\boldsymbol{B}_m} \right) \cdot \left( \sum_{n=1}^N \xi_n \operatorname{curl} (\operatorname{curl} \overline{\boldsymbol{B}_n}) \right) \,\dd x\,\dd t  \nonumber\\
    ={}& \varepsilon \sum_{|m-n|\leqslant 1} \int_0^T \psi(t) \int_\Omega \left[ \operatorname{curl}\overline{\boldsymbol{B}_n} \cdot \operatorname{curl}(\theta \xi_m \xi_n \overline{\boldsymbol{B}_m}) - \operatorname{div}\left( \operatorname{curl}\overline{\boldsymbol{B}_n} \times (\theta \xi_m \xi_n \overline{\boldsymbol{B}_m}) \right) \right] \,\dd x\,\dd t\nonumber\\
    ={}& \varepsilon \sum_{|m-n|\leqslant 1} \int_0^T \psi(t) \left[ \int_\Omega \operatorname{curl}\overline{\boldsymbol{B}_n} \cdot \operatorname{curl}(\theta \xi_m \xi_n \overline{\boldsymbol{B}_m}) \, dx - \int_{\partial\Omega} \left( \operatorname{curl}\overline{\boldsymbol{B}_n} \times (\theta \xi_m \xi_n \overline{\boldsymbol{B}_m}) \right) \cdot \boldsymbol{n} \, \dd s \right] \,\dd x\,\dd t \nonumber\\
    ={}& \varepsilon \sum_{|m-n|\leqslant 1} \int_0^T \psi(t) \int_\Omega \operatorname{curl}\overline{\boldsymbol{B}_n} \cdot \operatorname{curl}(\theta \xi_m \xi_n \overline{\boldsymbol{B}_m}) \,\dd x\,\dd t \nonumber\\
    ={}& \varepsilon \sum_{|m-n|\leqslant 1} \int_0^T \psi(t) \int_\Omega \operatorname{curl}\overline{\boldsymbol{B}_n} \cdot \left[ (\theta \xi_m \xi_n) \operatorname{curl}\overline{\boldsymbol{B}_m} + \nabla(\theta \xi_m \xi_n) \times \overline{\boldsymbol{B}_m} \right] \,\dd x\,\dd t\nonumber\\
    ={}& \varepsilon \sum_{|m-n|\leqslant 1} \int_0^T \psi(t) \int_\Omega \theta(x) \xi_m \xi_n \operatorname{curl}\overline{\boldsymbol{B}_n} \cdot \operatorname{curl}\overline{\boldsymbol{B}_m}\,\dd x\,\dd t  \nonumber\\
    &+ \varepsilon \sum_{|m-n|\leqslant 1} \int_0^T \psi(t) \int_\Omega \operatorname{curl}\overline{\boldsymbol{B}_n} \cdot (\nabla(\theta \xi_m \xi_n) \times \overline{\boldsymbol{B}_m}) \,\dd x\,\dd t \nonumber\\
    ={}& \varepsilon \sum_{|m-n|\leqslant 1} \int_0^T \psi(t) \left[ \int_\Omega \theta(x) \xi_m \xi_n \operatorname{curl}\overline{\boldsymbol{B}_n} \cdot \operatorname{curl}\overline{\boldsymbol{B}_m} \right. + \int_\Omega \theta(\nabla(\xi_m \xi_n) \times \overline{\boldsymbol{B}_m}) \cdot \operatorname{curl}\overline{\boldsymbol{B}_n} \nonumber\\
    &\qquad\qquad\quad \left. + \int_\Omega \xi_m \xi_n (\nabla\theta \times \overline{\boldsymbol{B}_m}) \cdot \operatorname{curl}\overline{\boldsymbol{B}_n} \right] \,\dd x\,\dd t \nonumber\\
    =:{}& R_1 + R_2 + R_3.
\end{align}

In the sequel, we shall estimate $R_i$, $i=1,2,3$ in order.
\vspace{1em}

\noindent
\underline{\textbf{Estimate for $R_1$}.} We consider two cases separately.

\smallskip
\noindent
\underline{Case~1: when $m=n$.} Since $V_n$ overlaps with the neighbouring patch $V_{n+1}$, we have
\begin{align}
    &\varepsilon\sum_{|m-n|=0} \int_0^T \psi \int_\Omega \theta \xi_m \xi_n {(\operatorname{curl}\overline{\mathbf B}_n)\cdot(\operatorname{curl}\overline{\mathbf B}_m)} \,\dd x\,\dd t\nonumber\\
    &\quad=\varepsilon \int_0^T \psi \int_\Omega \theta\xi^2_n |\operatorname{curl}\overline{\mathbf B}_n|^2\,\dd x\,\dd t\nonumber\\
    &\quad=\varepsilon \int_0^T \psi \left(\int_{V_n\cap V_{n+1}} + \int_{V_n\cap V_{n+1}^c}\right) \theta \xi_n^2 |\operatorname{curl}\overline{\mathbf B}_n|^2\,\dd x\,\dd t \nonumber\\
    &\quad\leq C\varepsilon \int_0^T \left( \|\operatorname{curl}\overline{\mathbf B}_n\|^2_{L^3(V_n\cap V_{n+1})} \|\xi_n^2\|_{L^3(V_n\cap V_{n+1})} + \|\operatorname{curl}\overline{\mathbf B}_n\|^2_{L^3(V_n\cap V_{n+1}^c)} \|\xi_n^2\|_{L^3(V_n\cap V_{n+1}^c)} \right) \,\dd t \nonumber\\
    &\quad\leq C\varepsilon \left( \varepsilon^{2\beta_{n+1}(\sigma_2-1)+\frac{1}{3}\beta_n} + \varepsilon^{2\beta_n(\sigma_2-1)+\frac{1}{3}\beta_{n-1}} \right) \int_0^T \|\mathbf B\|_{B_3^{\sigma_2, \infty}(\Omega^\varepsilon)}^2\,\dd t,
\end{align}
where the last inequality follows from the estimates in \eqref{2.18} and \eqref{3.3}.

\smallskip
\underline{Case~2: when $|m-n|=1$.} In this case, similarly, we have
\begin{align}
    &\varepsilon \sum_{|m-n|=1}\int_0^T \psi \int_\Omega \theta \xi_m \xi_n \operatorname{curl}\overline{\mathbf{B}}_n\cdot\operatorname{curl}\overline{\mathbf{B}}_m\,\dd x\,\dd t \nonumber \\
    ={}& \varepsilon \left[ \int_0^T \psi \int_\Omega \theta \xi_n \xi_{n+1} \operatorname{curl}\overline{\mathbf{B}}_n\cdot\operatorname{curl}\overline{\mathbf{B}}_{n+1}\,\dd x\,\dd t + \int_0^T \psi \int_\Omega \theta \xi_n \xi_{n-1} \operatorname{curl}\overline{\mathbf{B}}_n\cdot\operatorname{curl}\overline{\mathbf{B}}_{n-1}\,\dd x\,\dd t \right] \nonumber \\
    \leqslant{}& C\varepsilon \int_0^T \|\operatorname{curl}\overline{\mathbf{B}}_n\|_{L^3(V_n\cap V_{n+1})} \|\operatorname{curl}\overline{\mathbf{B}}_{n+1}\|_{L^3(V_n\cap V_{n+1})} \|\xi_n\xi_{n+1}\|_{L^3(V_n\cap V_{n+1})}\,\dd t \nonumber \\
    &+ C\varepsilon \int_0^T \|\operatorname{curl}\overline{\mathbf{B}}_n\|_{L^3(V_n\cap V_{n-1})} \|\operatorname{curl}\overline{\mathbf{B}}_{n-1}\|_{L^3(V_n\cap V_{n-1})} \|\xi_n\xi_{n-1}\|_{L^3(V_n\cap V_{n-1})}\,\dd t \nonumber \\
    \leqslant{}& C\varepsilon \left( \varepsilon^{2\beta_{n+1}(\sigma_2-1)+\frac{1}{3}\beta_n} + \varepsilon^{2\beta_n(\sigma_2-1)+\frac{1}{3}\beta_{n-1}} \right) \int_0^T \|\mathbf{B}\|_{B_3^{\sigma_2, \infty}(\Omega^\varepsilon)}^2\,\dd t.
\end{align}

By the construction of the sequence $\{\beta_n\}$ in (2.12), the exponent satisfies
\[
2\beta_n(\sigma_2-1)+\frac{1}{3}\beta_{n-1}+1>0.
\]
Hence, $R_1$ tends to zero as $\varepsilon \to 0$.

\vspace{1em}

\noindent
\underline{\textbf{Estimate for $R_2$}.} This term vanishes identically: in view of~\eqref{2.20} and~\eqref{3.2}, one has that
\begin{align}\label{R2}
    R_2 &= \varepsilon \sum_{|m-n|\leqslant 1} \int_0^T \psi \int_\Omega \theta \bigl(\nabla(\xi_m\xi_n)\times\overline{\mathbf B}_m\bigr)\cdot\operatorname{curl}\overline{\mathbf B}_n\,\dd x\,\dd t\nonumber \\
    &= \sum_{n=1}^N \varepsilon \int_0^T \psi \left(\int_{V_n\cap V_{n-1}} + \int_{V_n\cap V_{n+1}}\right) \theta\bigl(\nabla(\xi_n^2)\times\overline{\mathbf B}_n\bigr)\cdot\operatorname{curl}\overline{\mathbf B}_n\,\dd x\,\dd t \nonumber \\
    &\quad + \sum_{n=1}^{N-1} \varepsilon \int_0^T \psi \int_{V_n\cap V_{n+1}} \theta\bigl(\nabla(\xi_n\xi_{n+1})\times\overline{\mathbf B}_n\bigr)\cdot\operatorname{curl}\overline{\mathbf B}_{n+1} + \theta\bigl(\nabla(\xi_n\xi_{n+1})\times\overline{\mathbf B}_{n+1}\bigr)\cdot\operatorname{curl}\overline{\mathbf B}_n\,\dd x\,\dd t \nonumber \\
    &= \sum_{n=1}^N \varepsilon \int_0^T \psi \left( \int_{V_{n+1}\cap V_n} \theta\bigl(\nabla(\xi_{n+1}^2)\times\overline{\mathbf B}_{n+1}\bigr)\cdot\operatorname{curl}\overline{\mathbf B}_{n+1}\,\dd x + \int_{V_n\cap V_{n+1}} \theta\bigl(\nabla(\xi_n^2)\times\overline{\mathbf B}_n\bigr)\cdot\operatorname{curl}\overline{\mathbf B}_{n+1}\,\dd x \right) \,\dd t \nonumber \\
    &\quad + \sum_{n=1}^{N-1} \varepsilon \int_0^T \psi \int_{V_n\cap V_{n+1}} \theta\bigl(\nabla(2\xi_n\xi_{n+1})\times\overline{\mathbf B}_{n+1}\bigr)\cdot\operatorname{curl}\overline{\mathbf B}_{n+1}\,\dd x\,\dd t \nonumber \\
    &= \sum_{n=1}^{N-1} \varepsilon \int_0^T \psi \int_{V_n\cap V_{n+1}} \theta\bigl(\nabla(\xi_n+\xi_{n+1})^2\times\overline{\mathbf B}_{n+1}\bigr)\cdot\operatorname{curl}\overline{\mathbf B}_{n+1}\,\dd x\,\dd t = 0.
\end{align}

\vspace{1em}

\noindent
\underline{\textbf{Estimate for $R_3$}.} For the boundary term $R_3$, {recall from Section 3.8 that} $\nabla\theta$ is supported in the boundary layer $\Gamma_{4\varepsilon}$ and that $V_N$ contains the entire region where $\nabla\theta\neq 0$, so we restrict all subsequent integrals to $V_N$. Thus, by applying the trace inequality, we obtain that
\begin{align}
R_3 =& \varepsilon \sum_{|m-n|\leqslant 1} \int_0^T \psi \int_\Omega \xi_m \xi_n (\nabla\theta \times \overline{\boldsymbol{B}}_m)\cdot \operatorname{curl}\overline{\boldsymbol{B}}_n\,\dd x\,\dd t \nonumber \\
\leqslant{}& \varepsilon \int_0^T \psi \int_{\Omega^{2\varepsilon}\cap \Gamma_{4\varepsilon}} |\operatorname{curl}\overline{\boldsymbol{B}}_N|\,|\nabla\theta|\,\xi_N^2\,|\overline{\boldsymbol{B}}_N - \boldsymbol{B} + \boldsymbol{B}|\,\dd x\,\dd t  \nonumber \\
    \leqslant{}& \varepsilon \int_0^T \psi \int_{\Omega^{2\varepsilon}\cap \Gamma_{4\varepsilon}} |\operatorname{curl}\overline{\boldsymbol{B}}_N|\,|\nabla\theta|\,\xi_N^2\,|\overline{\boldsymbol{B}}_N - \boldsymbol{B}|\,\dd x\,\dd t  + \varepsilon \int_0^T \psi \int_{\Omega^{2\varepsilon}\cap \Gamma_{4\varepsilon}} |\operatorname{curl}\overline{\boldsymbol{B}}_N|\,\xi_N^2\,|\nabla\theta|\,|\boldsymbol{B}|\,\dd x\,\dd t  \nonumber \\
    \leqslant{}& \left(\varepsilon \int_0^T \psi \int_{\Omega^{2\varepsilon}\cap \Gamma_{4\varepsilon}} |\operatorname{curl}\overline{\boldsymbol{B}}_N|^2\,\dd x\,\dd t \right)^{\frac{1}{2}} \cdot \left(\varepsilon \int_0^T \psi \int_{\Omega^{2\varepsilon}\cap \Gamma_{4\varepsilon}} |\nabla\theta|^2 \xi_N^4 |\overline{\boldsymbol{B}}_N - \boldsymbol{B}|^2\,\dd x\,\dd t \right)^{\frac{1}{2}} \nonumber \\
    &+ \left(\varepsilon \int_0^T \psi \int_{\Omega^{2\varepsilon}\cap \Gamma_{4\varepsilon}} |\operatorname{curl}\overline{\boldsymbol{B}}_N|^2\,\dd x\,\dd t \right)^{\frac{1}{2}} \cdot \left(\varepsilon \int_0^T \psi \int_{\Omega^{2\varepsilon}\cap \Gamma_{4\varepsilon}} \xi_N^4 |\nabla\theta|^2 |\boldsymbol{B}|^2\,\dd x\,\dd t \right)^{\frac{1}{2}} \nonumber \\
    \leqslant{}& C \left(\varepsilon \int_0^T \int_\Omega |\operatorname{curl}\overline{\boldsymbol{B}}_N|^2\,\dd x\,\dd t \right)^{\frac{1}{2}} \cdot \left(\varepsilon \int_0^T \varepsilon^{-2} \|\xi_N^4\|_{L^3(\Omega^{2\varepsilon}\cap \Gamma_{4\varepsilon})} \|\overline{\boldsymbol{B}}_N - \boldsymbol{B}\|_{L^3(\Omega^{2\varepsilon}\cap \Gamma_{4\varepsilon})}^2\,dt\right)^{\frac{1}{2}} \nonumber \\
    &+ \left(\varepsilon \int_0^T \int_\Omega |\operatorname{curl}\overline{\boldsymbol{B}}_N|^2\,\dd x\,\dd t \right)^{\frac{1}{2}} \cdot \left(\varepsilon \int_0^T \|\nabla\theta \boldsymbol{B}\|_{L^2(\Omega^{2\varepsilon}\cap \Gamma_{4\varepsilon})}^2\,\dd t\right)^{\frac{1}{2}} \nonumber \\
    \leqslant{}& C \left(\varepsilon \int_0^T \int_\Omega |\operatorname{curl}\overline{\boldsymbol{B}}_N|^2\,\dd x\,\dd t \right)^{\frac{1}{2}} \cdot \left(\varepsilon \int_0^T \varepsilon^{2\sigma_2-\frac{5}{3}} \|\boldsymbol{B}\|_{B_3^{\sigma_2, \infty}(\Gamma_{4\varepsilon})}^2\,\dd t\right)^{\frac{1}{2}} \nonumber \\
    &+ \left(\varepsilon \int_0^T \int_\Omega |\operatorname{curl}\overline{\boldsymbol{B}}_N|^2\,\dd x\,\dd t \right)^{\frac{1}{2}} \cdot \left(\varepsilon \int_0^T \int_{\Gamma_{4\varepsilon}} |\nabla\overline{\boldsymbol{B}}_N|^2\,\dd t\right)^{\frac{1}{2}} \nonumber \\
    \leqslant{}& C \left[ \left(\varepsilon^{2\sigma_2-\frac{2}{3}}\right)^{\frac{1}{2}} + \left(\varepsilon \int_0^T \int_{\Gamma_{4\varepsilon}} |\nabla\overline{\boldsymbol{B}}_N|^2\,\dd x\,\dd t \right)^{\frac{1}{2}} \right] \longrightarrow 0.
\end{align}

Combining the estimates for $R_1, R_2$, and $R_3$, we complete the proof of Proposition~\ref{prop: 3.2}.
\end{proof}

It remains to deal with the most technical term, \textit{i.e.}, the last one, on the right-hand side of~\eqref{3.20}. To this end, we establish the following:

\begin{proposition}\label{prop: 3.3}
    Under the assumptions of Theorem~\ref{thm: main}, it holds that 
\begin{equation}\label{3.28}
    \lim_{\varepsilon \to 0} \left| \sum_{|m-n|\leqslant 1} \int_0^T \psi(t) \int_\Omega \theta(x) \xi_m \xi_n \left[ \overline{\boldsymbol{u}_m} \cdot \overline{(\operatorname{curl} \boldsymbol{B} \times \boldsymbol{B})_n} + \overline{\boldsymbol{B}_m} \cdot \operatorname{curl}\left(\overline{\boldsymbol{u} \times \boldsymbol{B}}\right)_n \right] \,\dd x\,\dd t \right| = 0.
\end{equation}
\end{proposition} 
\begin{proof}
We decompose the integral as follows:
\begin{align}\label{3.29}
    &\sum_{|m-n|\leqslant 1} \int_0^T \psi(t) \int_\Omega \theta(x) \xi_m \xi_n \left[ \overline{\boldsymbol{u}_m} \cdot \overline{(\operatorname{curl} \boldsymbol{B} \times \boldsymbol{B})_n} + \overline{\boldsymbol{B}_m} \cdot \operatorname{curl}\left(\overline{\boldsymbol{u} \times \boldsymbol{B}}\right)_n \right] \,\dd x\,\dd t \nonumber \\
    ={}& \sum_{|m-n|\leqslant 1} \int_0^T \psi(t) \int_\Omega \theta(x) \xi_m \xi_n \overline{\boldsymbol{u}_m} \cdot \left[ \operatorname{div}(\overline{\boldsymbol{B}_n} \otimes \overline{\boldsymbol{B}_n}) - \frac{1}{2} \nabla |\overline{\boldsymbol{B}_n}|^2 \right.\nonumber \\
    &\qquad\qquad\qquad\qquad\qquad \left. + \operatorname{div}\left(\overline{(\boldsymbol{B} \otimes \boldsymbol{B})_n} - \overline{\boldsymbol{B}_n} \otimes \overline{\boldsymbol{B}_n}\right) - \frac{1}{2} \nabla \left(\overline{(|\boldsymbol{B}|^2)_n} - |\overline{\boldsymbol{B}_n}|^2\right) \right] \,\dd x\,\dd t\nonumber \\
    &+ \sum_{|m-n|\leqslant 1} \int_0^T \psi(t) \int_\Omega \theta(x) \xi_m \xi_n \overline{\boldsymbol{B}_m} \cdot \left[ \operatorname{curl}(\overline{\boldsymbol{u}_n} \times \overline{\boldsymbol{B}_n}) + \operatorname{curl}\left(\overline{(\boldsymbol{u} \times \boldsymbol{B})_n} - \overline{\boldsymbol{u}_n} \times \overline{\boldsymbol{B}_n}\right) \right] \,\dd x\,\dd t\nonumber \\
    &=I_0+I_1+I_2+I_3,
    \end{align}
where
\begin{align*}
&I_0:= \sum_{|m-n|\leqslant 1} \int_0^T \psi \int_\Omega \theta \xi_m \xi_n \left[ \overline{\boldsymbol{u}_m} \cdot \left(\operatorname{div}(\overline{\boldsymbol{B}_n} \otimes \overline{\boldsymbol{B}_n}) - \frac{1}{2} \nabla |\overline{\boldsymbol{B}_n}|^2\right) + \overline{\boldsymbol{B}_m} \cdot \operatorname{curl}(\overline{\boldsymbol{u}_n} \times \overline{\boldsymbol{B}_n}) \right]\,\dd x\,\dd t,\\
    &I_1 := \sum_{|m-n|\leqslant 1} \int_0^T \psi \int_\Omega \theta \xi_m \xi_n \overline{\boldsymbol{u}_m} \cdot \operatorname{div}\left(\overline{(\boldsymbol{B} \otimes \boldsymbol{B})_n} - \overline{\boldsymbol{B}_n} \otimes \overline{\boldsymbol{B}_n}\right) \,\dd x\,\dd t,\\
    &I_2:= \frac{1}{2} \sum_{|m-n|\leqslant 1} \int_0^T \psi \int_\Omega \theta \xi_m \xi_n \overline{\boldsymbol{u}_m} \cdot \nabla \left( \overline{(|\boldsymbol{B}|^2)_n} - |\overline{\boldsymbol{B}_n}|^2 \right) \,\dd x\,\dd t,\\
    &{I_3:= \sum_{|m-n|\leqslant 1} \int_0^T \psi \int_\Omega \theta \xi_m \xi_n \overline{\boldsymbol{B}_m} \cdot \operatorname{curl}\left( \overline{(\boldsymbol{u} \times \boldsymbol{B})_n} - \overline{\boldsymbol{u}_n} \times \overline{\boldsymbol{B}_n} \right) \,\dd x\,\dd t.}
\end{align*}
    
\vspace{1em}

\noindent
\underline{\textbf{Estimate for $I_0$}.} This term vanishes identically by vector identities: 
\begin{align}
    I_0 :=&\sum_{|m-n|\leqslant 1} \int_0^T \psi \int_\Omega \theta \xi_m \xi_n \left[ \overline{\boldsymbol{u}_m} \cdot \left(\operatorname{div}(\overline{\boldsymbol{B}_n} \otimes \overline{\boldsymbol{B}_n}) - \frac{1}{2} \nabla |\overline{\boldsymbol{B}_n}|^2\right) + \overline{\boldsymbol{B}_m} \cdot \operatorname{curl}(\overline{\boldsymbol{u}_n} \times \overline{\boldsymbol{B}_n}) \right] \,\dd x\,\dd t \nonumber\\
    ={}& \sum_{|m-n|\leqslant 1} \int_0^T \psi \int_\Omega \theta \xi_m \xi_n \left[ \overline{\boldsymbol{u}_m} \cdot (\operatorname{curl}\overline{\boldsymbol{B}_n} \times \overline{\boldsymbol{B}_n}) + \overline{\boldsymbol{B}_m} \cdot \operatorname{curl}(\overline{\boldsymbol{u}_n} \times \overline{\boldsymbol{B}_n}) \right] \,\dd x\,\dd t \nonumber\\
    ={}& \sum_{|m-n|\leqslant 1} \int_0^T \psi \int_\Omega \theta \xi_m \xi_n \left[ \overline{\boldsymbol{u}_{\tilde{n}}} \cdot (\operatorname{curl}\overline{\boldsymbol{B}_{\tilde{n}}} \times \overline{\boldsymbol{B}_{\tilde{n}}}) + \operatorname{div}\left( ( \overline{\boldsymbol{u}_{\tilde{n}}} \times \overline{\boldsymbol{B}_{\tilde{n}})} \times \overline{\boldsymbol{B}_n} \right) - \overline{\boldsymbol{u}_{\tilde{n}}} \cdot (\operatorname{curl}\overline{\boldsymbol{B}_{\tilde{n}}} \times \overline{\boldsymbol{B}_{\tilde{n}}}) \right] \,\dd x\,\dd t \nonumber\\
    ={}&0. \label{3.30}
\end{align}
It remains to estimate the error and commutator terms $I_i$ $(i=1,2,3)$:

\vspace{1em}

\noindent
\underline{\textbf{Estimate for $I_1$}.} By integration by parts, we further decompose $I_1$ into three terms:
\begin{align}
    I_1 :=& \sum_{|m-n|\leqslant 1} \int_0^T \psi \int_\Omega \theta \xi_m \xi_n \overline{\boldsymbol{u}_m} \cdot \operatorname{div}\left(\overline{(\boldsymbol{B} \otimes \boldsymbol{B})_n} - \overline{\boldsymbol{B}_n} \otimes \overline{\boldsymbol{B}_n}\right) \,\dd x\,\dd t \nonumber\\
    ={}& -\sum_{|m-n|\leqslant 1} \int_0^T \psi \int_\Omega \nabla \theta \cdot (\xi_m \xi_n \overline{\boldsymbol{u}_m}) \left( \overline{(\boldsymbol{B} \otimes \boldsymbol{B})_n} - \overline{\boldsymbol{B}_n} \otimes \overline{\boldsymbol{B}_n} \right) \,\dd x\,\dd t \nonumber\\
    &- \sum_{|m-n|\leqslant 1} \int_0^T \psi \int_\Omega \theta \nabla (\xi_m \xi_n) \cdot \overline{\boldsymbol{u}_m} \left( \overline{(\boldsymbol{B} \otimes \boldsymbol{B})_n} - \overline{\boldsymbol{B}_n} \otimes \overline{\boldsymbol{B}_n} \right)\,\dd x\,\dd t\nonumber\\
    &- \sum_{|m-n|\leqslant 1} \int_0^T \psi \int_\Omega \theta \xi_m \xi_n \nabla \overline{\boldsymbol{u}_m} \cdot \left(\overline{(\boldsymbol{B} \otimes \boldsymbol{B})_n} - \overline{\boldsymbol{B}_n} \otimes \overline{\boldsymbol{B}_n}\right)\,\dd x\,\dd t =: I_{11}+I_{12}+I_{13}. \label{3.31}
\end{align}
\begin{itemize}
    \item 
For $I_{11}$, due to the presence of $\nabla \theta$, the nonzero contribution to the spatial integral comes only from the boundary layer $\Gamma_{4\e}$. Applying the Cauchy--Schwarz inequality together with \eqref{2.4} and \eqref{2.8}, we obtain that
\begin{align}
    |I_{11}|=& \left| \sum_{|m-n|\leqslant 1} \int_0^T \psi \int_\Omega \nabla \theta \cdot (\xi_m \xi_n \overline{\boldsymbol{u}_m}) \left( \overline{(\boldsymbol{B} \otimes \boldsymbol{B})_n} - \overline{\boldsymbol{B}_n} \otimes \overline{\boldsymbol{B}_n} \right) \,\dd x\,\dd t \right| \nonumber\\
    \leqslant{}& C \left( \int_0^T \int_{\Gamma_{4\varepsilon}} |\overline{\boldsymbol{B}_N}|^4 \,\dd x\,\dd t\right)^{\frac{1}{2}} \left( \int_0^T \int_{\Gamma_{4\varepsilon} \cap \Omega^\varepsilon} |\nabla \theta (\overline{\boldsymbol{u}_N} - \boldsymbol{u})|^2 \,\dd x\,\dd t\right)^{\frac{1}{2}} \nonumber\\
    \leqslant{}& C \left( \varepsilon \int_0^T \|\boldsymbol{B}\|_{L^4(\Gamma_{4\varepsilon})}^4 \,\dd t \right)^{\frac{1}{2}} \left( \int_0^T \varepsilon^{-\frac{5}{3}} \|\overline{\boldsymbol{u}_N} - \boldsymbol{u}\|_{L^3}^2 \,\dd t + \int_0^T \int_{\Gamma_{4\varepsilon}} |\nabla \theta \boldsymbol{u}|^2 \,\dd x\,\dd t\right)^{\frac{1}{2}} \nonumber\\
    \leqslant{}& C \varepsilon^{\frac{1}{2}} \left( \int_0^T \varepsilon^{\left(2\sigma_1 - \frac{5}{3}\right)} \|\boldsymbol{u}\|_{B_{3, \infty}^{\sigma_1}(\Omega^\varepsilon)}^2 \,\dd t + \int_0^T \int_{\Gamma_{4\varepsilon}} |\nabla \boldsymbol{u}|^2 \,\dd x\,\dd t \right)^{\frac{1}{2}} \nonumber\\
    \leqslant{}& C \left( \varepsilon^{\left(2\sigma_1 - \frac{2}{3}\right)} + \varepsilon \int_0^T \int_{\Gamma_{4\varepsilon}} |\nabla \boldsymbol{u}|^2 \,\dd x\,\dd t \right)^{\frac{1}{2}} \longrightarrow 0 \quad \text{as } \varepsilon \longrightarrow 0. \label{3.32}
\end{align}

\item 
The argument for  $I_{12}$ is parallel to that for $R_2$ above; see \eqref{R2}. Indeed,
\begin{align}
   I_{12}=
&\sum_{|m-n|\leqslant 1}
\int_0^T \psi \int_\Omega
\theta \nabla(\xi_m\xi_n)\cdot\overline{\boldsymbol{u}}_m
\left(
\overline{(\boldsymbol{B}\otimes\boldsymbol{B})}_n
-\overline{\boldsymbol{B}}_n\otimes\overline{\boldsymbol{B}}_n
\right)\,\dd x\,\dd t \nonumber
\\
&=
\sum_{n=1}^N
\int_0^T \psi
\left(
\int_{V_n\cap V_{n-1}}
+\int_{V_n\cap V_{n+1}}
\right)
\theta \nabla(\xi_n^2)\cdot\overline{\boldsymbol{u}}_n
\left(
\overline{(\boldsymbol{B}\otimes\boldsymbol{B})}_n
-\overline{\boldsymbol{B}}_n\otimes\overline{\boldsymbol{B}}_n
\right)\,\dd x\,\dd t \nonumber
\\
&\quad+
\sum_{n=1}^{N-1}
\int_0^T \psi
\int_{V_n\cap V_{n+1}}
\theta \nabla(\xi_n\xi_{n+1})\cdot\overline{\boldsymbol{u}}_n
\left(
\overline{(\boldsymbol{B}\otimes\boldsymbol{B})}_{n+1}
-\overline{\boldsymbol{B}}_{n+1}\otimes
\overline{\boldsymbol{B}}_{n+1}
\right) \nonumber
\\
&\qquad\qquad\qquad+
\theta \nabla(\xi_n\xi_{n+1})\cdot
\overline{\boldsymbol{u}}_{n+1}
\left(
\overline{(\boldsymbol{B}\otimes\boldsymbol{B})}_n
-\overline{\boldsymbol{B}}_n\otimes
\overline{\boldsymbol{B}}_n
\right)\,\dd x\,\dd t \nonumber
\\
&=
\sum_{n=1}^{N-1}
\int_0^T \psi
\int_{V_n\cap V_{n+1}}
\theta
\nabla(\xi_n^2+2\xi_n\xi_{n+1}+\xi_{n+1}^2)
\cdot\overline{\boldsymbol{u}}_{n+1} \left(
\overline{(\boldsymbol{B}\otimes\boldsymbol{B})}_{n+1}
-\overline{\boldsymbol{B}}_{n+1}\otimes
\overline{\boldsymbol{B}}_{n+1}
\right)\,\dd x\,\dd t \nonumber
\\
&=
\sum_{n=1}^{N-1}
\int_0^T \psi
\int_{V_n\cap V_{n+1}}
\theta
\nabla(\xi_n+\xi_{n+1})^2
\cdot\overline{\boldsymbol{u}}_{n+1}
\cdot
\left(
\overline{(\boldsymbol{B}\otimes\boldsymbol{B})}_{n+1}
-\overline{\boldsymbol{B}}_{n+1}\otimes
\overline{\boldsymbol{B}}_{n+1}
\right)\,\dd x\,\dd t\nonumber\\
&=0. \label{3.33}
\end{align}
In the above, the third equality follows from the construction~\eqref{3.2} on the overlapping region $V_n \cap V_{n+1}$, and the last equality follows from  $\nabla(\xi_n + \xi_{n+1})^2 = 0$.

\item 
For $I_{13}$, we consider two cases separately:

\smallskip
\noindent
\underline{Case~1: $m=n$.} From H\"older's inequality, we infer that
\begin{align}
    &\sum_{n=1}^{N-1}\int_0^T \psi \int_\Omega \theta \xi_n^2 \nabla \overline{\boldsymbol{u}_n} : \left( \overline{(\boldsymbol{B} \otimes \boldsymbol{B})_n} - \overline{\boldsymbol{B}_n} \otimes \overline{\boldsymbol{B}_n} \right)\,\dd x\,\dd t\nonumber\\
    \leqslant{}& C\sum_{n=1}^{N-1} \int_0^T \left\| \overline{(\boldsymbol{B} \otimes \boldsymbol{B})_n} - \overline{\boldsymbol{B}_n} \otimes \overline{\boldsymbol{B}_n} \right\|_{L^\frac{3}{2}(V_n \cap V_{n+1})} \|\nabla \overline{\boldsymbol{u}_n}\|_{L^3(V_n \cap V_{n+1})} \,\dd t \nonumber\\
    &+ C \sum_{n=1}^{N-1} \int_0^T \left\| \overline{(\boldsymbol{B} \otimes \boldsymbol{B})_n} - \overline{\boldsymbol{B}_n} \otimes \overline{\boldsymbol{B}_n} \right\|_{L^\frac{3}{2}(V_n \cap V_{n+1}^c)} \|\overline{\boldsymbol{u}_n}\|_{L^3(V_n \cap V_{n+1}^c)} \,\dd t \nonumber\\
    \leqslant{}& C \sum_{n=1}^{N-1} \int_0^T \left( \varepsilon^{2\beta_{n+1}\sigma_2 + \beta_{n+1}(\sigma_1-1)} + \varepsilon^{2\beta_n\sigma_2 + \beta_n(\sigma_1-1)} \right) \|\boldsymbol{B}\|_{B_3^{\sigma_2, \infty}(\Omega^\varepsilon)}^2 \|\boldsymbol{u}\|_{B_3^{\sigma_1, \infty}(\Omega^\varepsilon)} \,\dd t \nonumber\\
    \leqslant{}& C \sum_{n=1}^{N-1}\varepsilon^{\beta_n(2\sigma_2 + \sigma_1 - 1)} \longrightarrow 0 \quad \text{as } \varepsilon \longrightarrow 0. \label{3.34}
\end{align}

\smallskip
\noindent
\underline{Case~2: $|m-n|=1$.} In this case, we have that 
\begin{align}
&\left|
\sum_{|m-n|=1}
\int_0^T\int_\Omega
\theta\xi_m\xi_n
\nabla\overline{\boldsymbol{u}_m}:
\left(
\overline{(\boldsymbol{B}\otimes\boldsymbol{B})_n}
-
\overline{\boldsymbol{B}_n}\otimes
\overline{\boldsymbol{B}_n}
\right)\,\dd x\,\dd t
\right| \nonumber
\\
&\leqslant
C\sum_{|m-n|=1}
\left\|\nabla\overline{\boldsymbol{u}_m}\right\|_
{L^3(\widetilde{\Omega}_T)}
\left\|
\overline{(\boldsymbol{B}\otimes\boldsymbol{B})_n}
-
\overline{\boldsymbol{B}_n}\otimes
\overline{\boldsymbol{B}_n}
\right\|_{L^{\frac 32}(\widetilde{\Omega}_T)} \nonumber
\\
&\leqslant
C\sum_{|m-n|=1}
\varepsilon^{\beta_{\tilde n}(\sigma_1-1)}
\left\|\boldsymbol{u}^\varepsilon\right\|_
{B_3^{\sigma_1,\infty}(\widetilde{\Omega}_T)}
\varepsilon^{2\beta_{\tilde n}\sigma_2}
\left\|\boldsymbol{B}^\varepsilon\right\|_
{B_3^{\sigma_2,\infty}(\widetilde{\Omega}_T)}^2 \nonumber
\\
&\leqslant
C\sum_{|m-n|=1}
\varepsilon^{\beta_{\tilde n}(\sigma_1-1)
+2\beta_{\tilde n}\sigma_2} \nonumber
\\
&=
C\sum_{|m-n|=1}
\varepsilon^{\beta_{\tilde n}(2\sigma_2+\sigma_1-1)}
\longrightarrow 0 \quad \text{as } \varepsilon \longrightarrow 0.\label{3.35}
\end{align}

\end{itemize}

Putting together the bounds in~\eqref{3.31}--\eqref{3.35}, we deduce that $I_{1} \to 0$ as $\e\to 0$ under the same assumptions as in Theorem~\ref{thm: main}.

\vspace{1em}

\noindent
\underline{\textbf{Estimate for $I_2$}.}  Using integration by parts, we split $I_{2}$ into three terms:
\begin{align}
    I_{2} :=& \frac{1}{2} \sum_{|m-n|\leqslant 1} \int_0^T \psi \int_\Omega \theta \xi_m \xi_n \overline{\boldsymbol{u}_m} \cdot \nabla \left( \overline{(|\boldsymbol{B}|^2)_n} - |\overline{\boldsymbol{B}_n}|^2 \right) \,\dd x\,\dd t \nonumber\\
    ={}& -\frac{1}{2} \sum_{|m-n|\leqslant 1} \int_0^T \psi \int_\Omega \nabla\theta \cdot (\xi_m \xi_n \overline{\boldsymbol{u}_m}) \left( \overline{(|\boldsymbol{B}|^2)_n} - |\overline{\boldsymbol{B}_n}|^2 \right) \,\dd x\,\dd t \nonumber\\
    &-\frac{1}{2} \sum_{|m-n|\leqslant 1} \int_0^T \psi \int_\Omega \theta \nabla(\xi_m \xi_n) \cdot \overline{\boldsymbol{u}_m} \left( \overline{(|\boldsymbol{B}|^2)_n} - |\overline{\boldsymbol{B}_n}|^2 \right) \,\dd x\,\dd t \nonumber\\
    &-\frac{1}{2} \sum_{|m-n|\leqslant 1} \int_0^T \psi \int_\Omega \theta \xi_m \xi_n \operatorname{div}(\overline{\boldsymbol{u}_m}) \left( \overline{(|\boldsymbol{B}|^2)_n} - |\overline{\boldsymbol{B}_n}|^2 \right) \,\dd x\,\dd t=:I_{21}+I_{22}+I_{23}.\label{3.36}
\end{align}
By an analogous argument to that for $I_{12}$, we have $I_{22} \equiv 0$ thanks to the properties of the partition of unity, especially~\eqref{3.2}. To bound $I_{21}$, we note that $|\nabla\theta(x)| \lesssim \frac{1}{{\rm dist}(x,\p\Omega)}$ and estimate $$
\left\|\frac{\overline{\boldsymbol{u}_N}}{{\rm dist}(\cdot,\partial\Omega)}\right\|_{L^3{(\Gamma_{4\varepsilon})}} \lesssim \left\|\nabla\overline{\boldsymbol{u}_N}\right\|_{L^3{(\Gamma_{4\varepsilon})}}
$$
by Lemma~\ref{lem: hardy}. Meanwhile, for $I_{23}$, we observe that
$$
\left\| \overline{(|\boldsymbol{B}|^2)_n} - |\overline{\boldsymbol{B}_n}|^2 \right\|_{L^{\frac{3}{2}}{(V_n)}} \lesssim \varepsilon^{2\beta_n \sigma_2} \|\boldsymbol{B}\|_{B_3^{\sigma_2, \infty}(\Omega^\varepsilon)}^2.
$$
We may thus bound
\begin{align}
    |I_2| &\leqslant C \sum_{|m-n|\leqslant 1} \int_0^T \|\nabla \overline{\boldsymbol{u}_m}\|_{L^3{(V_m\cap V_n)}}  \left\| \overline{(|\boldsymbol{B}|^2)_n} - |\overline{\boldsymbol{B}_n}|^2 \right\|_{L^{\frac{3}{2}}{(V_m\cap V_n)}}\, \dd t \nonumber\\
    &\leqslant C \sum_{|m-n|\leqslant 1} \int_0^T \left( \varepsilon^{\beta_m(\sigma_1 - 1)} \|\boldsymbol{u}\|_{B_3^{\sigma_1, \infty}} \right) \left( \varepsilon^{2\beta_n \sigma_2} \|\boldsymbol{B}\|_{B_3^{\sigma_2, \infty}}^2 \right) \,\dd t \nonumber\\
    &\leqslant C \sum_{n=1}^N \varepsilon^{\beta_n(2\sigma_2 + \sigma_1 - 1)} \longrightarrow 0 \quad \text{as } \varepsilon \longrightarrow 0\text{ for $2\sigma_2 + \sigma_1 > 1$.} \label{3.37}
\end{align}
\smallskip


\noindent
\underline{\textbf{Estimate for $I_3$}.} This is the magnetic induction commutator term. Via integration by parts, we again split it into three terms:
\begin{align}
    I_3 :=& \sum_{|m-n|\leqslant 1} \int_0^T \psi \int_\Omega \theta \xi_m \xi_n \overline{\boldsymbol{B}_m} \cdot \operatorname{curl}\left( \overline{(\boldsymbol{u} \times \boldsymbol{B})_n} - \overline{\boldsymbol{u}_n} \times \overline{\boldsymbol{B}_n} \right) \,\dd x\,\dd t \nonumber\\
    ={}& \sum_{|m-n|\leqslant 1} \int_0^T \psi \int_\Omega \left[ \theta \xi_m \xi_n (\operatorname{curl}\overline{\boldsymbol{B}_m}) + \nabla(\theta \xi_m \xi_n) \times \overline{\boldsymbol{B}_m} \right] \cdot \left( \overline{(\boldsymbol{u} \times \boldsymbol{B})_n} - \overline{\boldsymbol{u}_n} \times \overline{\boldsymbol{B}_n} \right)  \,\dd x\,\dd t \nonumber\\
    ={}&\sum_{|m-n|\leqslant 1} \int_0^T \psi \int_\Omega \theta \xi_m \xi_n (\operatorname{curl}\overline{\boldsymbol{B}_m})\cdot \left( \overline{(\boldsymbol{u} \times \boldsymbol{B})_n} - \overline{\boldsymbol{u}_n} \times \overline{\boldsymbol{B}_n} \right) \,\dd x\,\dd t\nonumber\\
    &+\sum_{|m-n|\leqslant 1} \int_0^T \psi \int_\Omega \xi_m \xi_n
    (\nabla\theta  \times \overline{\boldsymbol{B}_m})\cdot \left( \overline{(\boldsymbol{u} \times \boldsymbol{B})_n} - \overline{\boldsymbol{u}_n} \times \overline{\boldsymbol{B}_n} \right)  \,\dd x\,\dd t \nonumber\\
    &+\sum_{|m-n|\leqslant 1} \int_0^T \psi \int_\Omega \theta \left(\nabla(\xi_m \xi_n) \times \overline{\boldsymbol{B}_m}\right) \cdot \left( \overline{(\boldsymbol{u} \times \boldsymbol{B})_n} - \overline{\boldsymbol{u}_n} \times \overline{\boldsymbol{B}_n} \right)  \,\dd x\,\dd t \nonumber\\
    =: {}&I_{31}+I_{32}+I_{33}. \label{3.38}
\end{align}

\begin{itemize}
    \item 
For $I_{31}$, from H\"{o}lder's inequality and the commutator estimate in Lemma~\ref{lem:commutator} we infer that
\begin{align}
    |I_{31}| ={}&\left|\sum_{|m-n|\leqslant 1} \int_0^T \psi \int_\Omega \theta \xi_m \xi_n (\operatorname{curl}\overline{\boldsymbol{B}_m}) \cdot \left( \overline{(\boldsymbol{u} \times \boldsymbol{B})_n} - \overline{\boldsymbol{u}_n} \times \overline{\boldsymbol{B}_n} \right) \,\dd x\,\dd t\right| \nonumber\\
    \leqslant{}& \sum_{|m-n|\leqslant 1} \int_0^T \psi \|\operatorname{curl}\overline{\boldsymbol{B}_m}\|_{L^3{(V_m\cap V_n)}} \left\| \overline{(\boldsymbol{u} \times \boldsymbol{B})_n} - \overline{\boldsymbol{u}_n} \times \overline{\boldsymbol{B}_n} \right\|_{L^{\frac{3}{2}}{(V_m\cap V_n)}} \,\dd t \nonumber\\
    \leqslant{}& C \sum_{|m-n|\leqslant 1} \int_0^T \left( \varepsilon^{\beta_m(\sigma_2 - 1)} \|\boldsymbol{B}\|_{B_3^{\sigma_2, \infty}{(\Omega^\varepsilon)}} \right) \left( \varepsilon^{\beta_n(\sigma_1 + \sigma_2)} \|\boldsymbol{u}\|_{B_3^{\sigma_1, \infty}{(\Omega^\varepsilon)}} \|\boldsymbol{B}\|_{B_3^{\sigma_2, \infty}{(\Omega^\varepsilon)}} \right) \,\dd t \nonumber\\
    \leqslant{}& C \sum_{|m-n|\leqslant 1} \varepsilon^{\beta_{\tilde{n}}(2\sigma_2 + \sigma_1 - 1)} \longrightarrow 0 \quad \text{as } \varepsilon \longrightarrow 0,  \label{3.39}
\end{align}
thanks to $2\sigma_2 + \sigma_1 > 1$.

\item 
For $I_{32}$, we adapt the estimates for $I_{22}$ to obtain that 
\begin{align}
    |I_{32}| ={}&\left|\sum_{|m-n|\leqslant 1} \int_0^T \psi \int_\Omega \xi_m \xi_n (\nabla\theta \times \overline{\boldsymbol{B}_m}) \cdot \left( \overline{(\boldsymbol{u} \times \boldsymbol{B})_n} - \overline{\boldsymbol{u}_n} \times \overline{\boldsymbol{B}_n} \right) \,\dd x\,\dd t \right| \nonumber\\
    \leqslant{}& C \int_0^T \psi \left\| \nabla\theta \times \overline{\boldsymbol{B}_N}\right\|_{L^3(\Gamma_{4\varepsilon})} \left\| \overline{(\boldsymbol{u} \times \boldsymbol{B})_N} - \overline{\boldsymbol{u}_N} \times \overline{\boldsymbol{B}_N} \right\|_{L^{\frac{3}{2}}{(\Gamma_{4\varepsilon})}} \,\dd t \nonumber\\
    \leqslant{}& C \int_0^T \psi \left\| \frac{\overline{\boldsymbol{B}_N}}{\operatorname{dist}(\cdot,\partial\Omega)} \right\|_{L^3(\Gamma_{4\varepsilon})} \left\| \overline{(\boldsymbol{u} \times \boldsymbol{B})_N} - \overline{\boldsymbol{u}_N} \times \overline{\boldsymbol{B}_N} \right\|_{L^{\frac{3}{2}}{(\Gamma_{4\varepsilon})}} \,\dd t \nonumber\\
    \leqslant{}& C \int_0^T \psi \|\nabla \overline{\boldsymbol{B}_N}\|_{L^3{(\Gamma_{4\varepsilon})}} \left\| \overline{(\boldsymbol{u} \times \boldsymbol{B})_N} - \overline{\boldsymbol{u}_N} \times \overline{\boldsymbol{B}_N} \right\|_{L^{\frac{3}{2}}{(\Gamma_{4\varepsilon})}} \,\dd t \nonumber\\
    \leqslant{}&  C \int_0^T \left( \varepsilon^{\sigma_2 - 1} \|\boldsymbol{B}\|_{B_3^{\sigma_2, \infty}{(\Omega^\varepsilon)}} \right) \left( \varepsilon^{\sigma_1 + \sigma_2} \|\boldsymbol{u}\|_{B_3^{\sigma_1, \infty}{(\Omega^\varepsilon)}} \|\boldsymbol{B}\|_{B_3^{\sigma_2, \infty}{(\Omega^\varepsilon)}} \right)\, \dd t \nonumber\\
    \leqslant{}& C \varepsilon^{2\sigma_2+\sigma_1 - 1}  \longrightarrow 0 \quad \text{as } \varepsilon \longrightarrow 0. \label{3.40}
\end{align}
We have applied H\"{o}lder's inequality, that $\nabla\theta$ is supported in the boundary layer $\Gamma_{4\e}$, Hardy's inequality (Lemma~\ref{lem: hardy}), the commutator estimate in Lemma~\ref{lem:commutator}, as well as $2\sigma_2+\sigma_1>1$.

\item 
Finally, by an analogous argument to ~\eqref{R2}, the properties of the  partition of unity $\{\xi_n\}$ lead to 
\begin{align}
   |I_3| &= \left|\sum_{|m-n|\leqslant 1} \int_0^T \psi \int_\Omega \theta \left(\nabla(\xi_m \xi_n) \times \overline{\boldsymbol{B}_m}\right) \cdot \left( \overline{(\boldsymbol{u} \times \boldsymbol{B})_n} - \overline{\boldsymbol{u}_n} \times \overline{\boldsymbol{B}_n} \right) \,\dd x\,\dd t\right|\nonumber\\
   &\qquad \qquad \longrightarrow 0 \quad \text{as } \varepsilon \longrightarrow 0. \label{3.41}
\end{align}
\end{itemize}

Combining the estimates in \eqref{3.29}--\eqref{3.41}, we complete the proof of Proposition~\ref{prop: 3.3}.  \end{proof}

With the preceding Propositions~\ref{prop: 3.1}--\ref{prop: 3.3} at hand, we deduce from the identity~\eqref{3.20} that 

\begin{align}
    \lim_{\varepsilon\to0} \sum_{|m-n|\le 1} \int_0^T \psi' \int_\Omega \theta \xi_m \xi_n \bigg( \overline{\rho}_{\tilde{n}} \frac{|\overline{\mathbf{u}}_{\tilde{n}}|^2}{2} + \frac{(\overline{\rho}_{\tilde{n}})^\gamma}{\gamma-1} + \frac{|\overline{\mathbf{B}}_{\tilde{n}}|^2}{2} \bigg) \,\dd x\,\dd t =0 \label{3.42}
\end{align}
The total energy of the MHD system at time $t$ is 
\[
\mathcal{E}^\varepsilon(t)
=\int_\Omega
\left(\frac{1}{2}\rho^\varepsilon |\mathbf{u}^\varepsilon|^2
+\frac{(\rho^\varepsilon)^\gamma}{\gamma-1}
+\frac{1}{2}|\mathbf{B}^\varepsilon|^2\right)\,\dd x,
\]
whose approximation by mollification discretises into the sum on the left-hand side of \eqref{3.42}. In view of the compressible MHD
system \eqref{1.2}, for arbitrary $\psi \in C^1_0(0,T)$ one has that
\begin{align}\label{eq:global_dissipation_bound}
    0
    &\le
    \varepsilon\int_0^T\psi\int_\Omega
    \left(
    |\nabla\mathbf{u}^\varepsilon|^2
    +
    |\nabla\mathbf{B}^\varepsilon|^2
    \right)\,\dd x\,\dd t
    \nonumber\\
    &\le
    \varepsilon\int_0^T\psi'\mathcal{E}^\varepsilon(t)\,\dd t
    +\int_0^T \psi(t) \int_\Omega \theta(x) \bigg(\sum_{m=1}^N \xi_m \overline{\mathbf{u}}_m\bigg) \cdot \left(\sum_{n=1}^N \xi_n \bigg(\partial_t (\overline{\rho\mathbf{u}})_n + \operatorname{div}\overline{(\rho\mathbf{u}\otimes \mathbf{u})}_n 
    +\nabla\overline{(\rho^\gamma)}_n\bigg)\right) \,\dd x\,\dd t \nonumber\\
    &+\int_0^T \psi(t) \int_\Omega \theta(x) \bigg(\sum_{m=1}^N \xi_m \overline{\mathbf{B}}_m\bigg) \cdot \left(\sum_{n=1}^N \xi_n \bigg(\partial_t \overline{\mathbf{B}}_n - \operatorname{curl}\overline{(\mathbf{u}\otimes \mathbf{B})}_n \bigg)\right) \,\dd x\,\dd t \nonumber\\
    &- \varepsilon \int_0^T \psi(t) \int_\Omega \theta(x) \bigg(\sum_{m=1}^N \xi_m \overline{\mathbf{u}}_m\bigg) \cdot \bigg(\sum_{n=1}^N \xi_n \Delta \overline{\mathbf{u}}_n\bigg) \,\dd x\,\dd t\nonumber\\
    &+\varepsilon \int_0^T \psi(t) \int_\Omega \theta(x) \bigg(\sum_{m=1}^N \xi_m \overline{\mathbf{B}}_m\bigg) \cdot \bigg(\sum_{n=1}^N \xi_n \operatorname{curl}(\operatorname{curl}\overline{\mathbf{B}}_n)\bigg) \,\dd x\,\dd t. 
\end{align}
By virtue of \eqref{3.10}, \eqref{3.20}, and Propositions~\ref{prop: 3.1}--\ref{prop: 3.3}, we conclude that the right-hand side of \eqref{eq:global_dissipation_bound} tends to zero when $\varepsilon \to 0$. This  completes the proof of the first part of Theorem~\ref{thm: main}.

\section{Convergence of density, momentum, and magnetic field}\label{sec: proof, convergence}

This section is devoted to the proof of the second part of Theorem~\ref{thm: main}; \textit{i.e.}, \eqref{1.9} reproduced below:
   \begin{equation*}
        \rho^\varepsilon \longrightarrow \rho^0, \quad \rho^\varepsilon\mathbf{u}^\varepsilon \longrightarrow \rho^0\mathbf{u}^0, \quad \mathbf{B}^\varepsilon \longrightarrow \mathbf{B}^0 \quad \text{in } L^3\left(0, T; L^p_{\rm loc}(\Omega)\right),\,1<p\leq \frac{9}{2}.
    \end{equation*} 
   We shall extract weakly converging subsequences and utilise the Aubin--Lions compactness lemma to pass to the limit in nonlinear terms.

\begin{proof}
    
 First, by the uniform energy bounds and the Besov regularity assumption~\eqref{1.5}, there exists $(\rho^0, \mathbf{u}^0, \mathbf{B}^0)$ such that for any interior domain $\widetilde{\Omega}_T \subset\subset \Omega \times (0,T)$, up to a subsequence (unrelabelled),
\begin{equation}\label{eq:weak_conv}
\begin{aligned}
    \rho^\varepsilon \rightharpoonup \rho^0 &\quad \text{weakly in } B_3^{\sigma_1,\infty}(\widetilde{\Omega}_T) \cap L^\infty(\Omega_T), \\
    \mathbf{u}^\varepsilon \rightharpoonup \mathbf{u}^0, \quad \rho^\varepsilon \mathbf{u}^\varepsilon \rightharpoonup \mathbf{M} &\quad \text{weakly in } B_3^{\sigma_1,\infty}(\widetilde{\Omega}_T), \\
    \mathbf{B}^\varepsilon \rightharpoonup \mathbf{B}^0 &\quad \text{weakly in } B_3^{\sigma_2,\infty}(\widetilde{\Omega}_T).
\end{aligned}
\end{equation}
From the global energy inequality, we have the dissipation bounds:
\begin{equation}\label{eq:dissipation_bound}
    \sqrt{\varepsilon}\nabla\mathbf{u}^\varepsilon \in L^2(0,T; L^2(\Omega)), \quad \sqrt{\varepsilon}\nabla\mathbf{B}^\varepsilon \in L^2(0,T; L^2(\Omega)).
\end{equation}
This together with \eqref{1.4}, \eqref{1.5} yields that for any spatial interior domain $\widetilde{\Omega} \subset\subset \Omega$ and time interval $I \subset\subset [0,T]$, one has that
\begin{equation}\label{eq:time_derivatives}
    \partial_t \rho^\varepsilon, \ \partial_t (\rho^\varepsilon \mathbf{u}^\varepsilon), \ \partial_t \mathbf{B}^\varepsilon \in L^2\left(I; W^{-1, \frac{3}{2}}(\widetilde{\Omega})\right).
\end{equation}
Recall the following continuous embeddings:
\begin{equation*}
    B_3^{\sigma_i,\infty}\left(I \times \widetilde{\Omega}\right) \hookrightarrow L^3\left(I; B_3^{\sigma_i,\infty}(\widetilde{\Omega})\right), \quad i=1,2,
\end{equation*}
as well as the compact embeddings:
\begin{equation}\label{eq:besov_embedding}
    B_3^{\sigma_i,\infty}(\widetilde{\Omega}) \hookrightarrow W^{\frac{1}{3}, \frac{9-\delta}{4-3\sigma_i}}(\widetilde{\Omega}) \hookrightarrow L^{\frac{9}{2}}(\widetilde{\Omega}),
\end{equation}
where $0 < \delta < \min\{9(\sigma_1 - \frac{1}{3}), 9(\sigma_2 - \frac{1}{3})\}$ and $\sigma_1, \sigma_2 > \frac{1}{3}$.

Thus, the assumption~\eqref{1.5}, the above regularity result in~\eqref{eq:time_derivatives}, and the compact embeddings in~\eqref{eq:besov_embedding}, we may apply the Aubin--Lions lemma to obtain that 
\begin{equation}\label{eq:strong_rho_B}
    \rho^\varepsilon \longrightarrow \rho^0, \quad \mathbf{B}^\varepsilon \longrightarrow \mathbf{B}^0 \quad \text{strongly in } L^3\left(I; L^p(\widetilde{\Omega})\right)\,\text{for any $p \in \left(1, \frac{9}{2}\right]$}.
\end{equation}
This combined with the weak convergence of $\mathbf{u}^\varepsilon$ in~\eqref{eq:weak_conv} implies $\mathbf{M} = \rho^0\mathbf{u}^0$ by the uniqueness of limits.  Similarly to \eqref{eq:strong_rho_B}, we deduce that 
\begin{equation}\label{eq:strong_momentum}
    \rho^\varepsilon \mathbf{u}^\varepsilon \longrightarrow \rho^0\mathbf{u}^0 \quad \text{strongly in } L^3\left(I; L^p(\widetilde{\Omega})\right).
\end{equation}
The convergence results in  \eqref{eq:strong_rho_B} and \eqref{eq:strong_momentum} allow us to pass to the limit in nonlinear terms: 
\begin{equation}\label{eq:nonlinear_limits}
\begin{aligned}
    \rho^\varepsilon \mathbf{u}^\varepsilon \otimes \mathbf{u}^\varepsilon &\longrightarrow \rho^0\mathbf{u}^0 \otimes \mathbf{u}^0, \\
    (\rho^\varepsilon)^\gamma &\longrightarrow (\rho^0)^\gamma, \\
    \mathbf{B}^\varepsilon \otimes \mathbf{B}^\varepsilon &\longrightarrow \mathbf{B}^0 \otimes \mathbf{B}^0, \\
    \mathbf{u}^\varepsilon \times \mathbf{B}^\varepsilon &\longrightarrow \mathbf{u}^0 \times \mathbf{B}^0\qquad \text{in $\mathcal{D}'\left(\widetilde{\Omega}_T\right)$}.
\end{aligned}
\end{equation}

Meanwhile, assuming the Kato-type condition~\eqref{1.7}, it is straightforward to see by Cauchy--Schwarz that for any $\boldsymbol{\psi} \in C_c^1(\Omega_T; \mathbb{R}^3)$, we have as $\e\to 0$ that
\begin{align}\label{eq:dissipation_vanish}
   & \e\left| \int_0^T \int_\Omega \nabla\mathbf{u}^\varepsilon : \nabla \boldsymbol{\psi}
    +\nabla\mathbf{B}^\varepsilon : \nabla \boldsymbol{\psi}\,\dd x\,\dd t \right| \nonumber\\
    &\qquad \leqslant C \varepsilon^{\frac{1}{2}} \left( \varepsilon \int_0^T \int_\Omega |\nabla \mathbf{u}^\varepsilon|^2
    +|\nabla \mathbf{B}^\varepsilon|^2\,\dd x\,\dd t \right)^{\frac{1}{2}} \longrightarrow 0.
\end{align}

Having established \eqref{eq:nonlinear_limits} and \eqref{eq:dissipation_vanish}, we can test the compressible MHD equations \eqref{1.1} against smooth functions compactly supported in $\Omega_T$, and take the limit $\varepsilon \rightarrow0$ to conclude that $(\rho^0, \mathbf{u}^0, \mathbf{B}^0)$ satisfy the ideal compressible MHD equation~\eqref{1.3} in the sense of distributions. 
\end{proof}

The proof of Theorem~\ref{thm: main} is thus complete.

\bigskip
\noindent
{\bf Acknowledgement}. SL is supported by NSFC Projects 12331008 $\&$ 12411530065, the Young Elite Scientists Sponsorship Program by CAST 2023QNRC001, National Key Research $\&$ Development Programs 2023YFA1010900 and 2024YFA1014900, Shanghai Rising-Star Program 24QA2703600, Shanghai Qi-Guang Scholarship, and Shanghai Frontiers Science Center of Modern Analysis.  TT is supported by NSFC No. 12371246, Qing Lan Project of Jiangsu Province.

\vspace{0.2cm}	
	
	\noindent
	{\bf Statement of competing interests}. 
	The authors declare that there is no conflict of interest.

	\noindent
	{\bf Statement of data availability}.
	Our manuscript has no associated data.

\noindent
{\bf AI Statement}. The authors use ChatGPT 6 models to produce Figure~1. All the mathematical derivations, proofs, and writing up of the manuscript have been worked out without the use of AI tools and remain solely the
responsibility of the authors.


\end{document}